

\input amstex
\documentstyle{amsppt}
\magnification=1200
\catcode`\@=11
\redefine\logo@{}
\catcode`\@=13

\define \bn{\Bbb N}
\define \bz{\Bbb Z}
\define \bq{\Bbb Q}
\define \br{\Bbb R}
\define \bc{\Bbb C}
\define \bh{\Bbb H_2}

\define \spi{Sp_4(\bz)}
\define \M{{\Cal M}}
\define\Ha{{\Cal H}}
\define\La{{\Cal L}}
\define\geg{{\goth g}}
\define\0o{{\overline 0}}
\define\1o{{\overline 1}}
\define\rk{\text{rk}~}
\define\gi{\Gamma_{\infty}}


\define\mult{\text{mult}}
\define\re{\text{re}}
\define\im{\text{im}}

\define\dtf{\Delta_5(Z)}
\define\dtp{[\Delta_5(Z)]_p}

\define\gm{\Gamma_1}
\define\gmt{\Gamma_t}
\define\slt{SL_2(\bz)}

\TagsOnRight

\document
\document

\topmatter
\title
The Igusa modular forms and ``the simplest'' Lorentzian
Kac--Moody algebras
\endtitle

\author
Valeri A. Gritsenko \footnote{Supported by
DFG grant Rus 17/108/95 and
Ruprecht--Karls--Universit\"at Heidelberg. \hfill\hfill} and
Viacheslav V. Nikulin \footnote{Partially supported by
Grant of Russian Fund of Fundamental Research and
Grant of AMS.
\hfill\hfill}
\endauthor

\address
St. Petersburg Department of Steklov Mathematical Institute,
\newline
${}\hskip 8pt $
Fontanka 27, 191011 St. Petersburg,  Russia
\endaddress
\email
gritsenk\@gauss.pdmi.ras.ru, valera\@mathi.uni-heidelberg.de
\endemail

\address
Steklov Mathematical Institute,
ul. Vavilova 42, Moscow 117966, GSP-1, Russia
\endaddress

\email
slava\@nikulin.mian.su
\endemail

\abstract
We find automorphic corrections for the Lorentzian
Kac--Moody algebras with the simplest
generalized Cartan matrices of  rank $3$
$$
A_{1,0}=
\pmatrix
\hphantom{-}{2}&\hphantom{-}{0}&{-1}&\cr
\hphantom{-}{0}&\hphantom{-}{2}&{-2}&\cr
{-1}&{-2}&\hphantom{-}{2}&\cr
\endpmatrix \ \ \ \text{and}\ \ \
A_{1,I}=
\pmatrix
\hphantom{-}{2}&{-2}&{-1}&\cr
{-2}&\hphantom{-}{2}&{-1}&\cr
{-1}&{-1}&\hphantom{-}{2}&\cr
\endpmatrix .
$$
For $A_{1,0}$ this correction,
which is a generalized Kac-Moody Lie superalgebra,
is given by the Igusa  $Sp_4(\bz)$-modular form $\chi_{35}(Z)$
of weight $35$,  and for $A_{1,I}$ by some Siegel modular  form
$\widetilde{\Delta}_{30}(Z)$ of weight $30$ with respect to
a $2$-congruence subgroup of $Sp_4(\bz)$.
We find the infinite product expansions for $\chi_{35}(Z)$
and $\widetilde{\Delta}_{30}(Z)$  and calculate multiplicities
of all roots of the corresponding generalized Lorentzian
Kac--Moody superalgebras. These multiplicities are
given by Fourier coefficients of some
Jacobi forms of weight $0$ and index $1$.

Our method of construction of $\chi_{35}(Z)$ and
$\widetilde{\Delta}_{30}(Z)$ naturally leads to the general
and direct construction by infinite product or sum expansions
of Siegel modular forms, whose divisors are
the Humbert surfaces with fixed discriminants. Existence of
these forms was proved by van der Geer in 1982 using some
geometrical considerations.

To show  the perspective for further study, we announce
a list of all hyperbolic symmetric
generalized Cartan matrices $A$ of rank 3 such that $A$ has
elliptic or parabolic type, $A$ has a lattice Weyl vector, and
$A$ contains the parabolic submatrix
$\widetilde{\Bbb A}_1$.
\endabstract

\rightheadtext
{Igusa modular  form and Kac--Moody algebras}
\leftheadtext{V. Gritsenko and  V. Nikulin}
\endtopmatter

\document
\head
0. Introduction
\endhead

The main starting point of this paper is to find
an automorphic correction of the Lorentzian Kac--Moody algebra
$\geg(A_{1,0})$ with the generalized Cartan matrix
$$
A_{1,0}=
\pmatrix
\hphantom{-}{2}&\hphantom{-}{0}&{-1}&\cr
\hphantom{-}{0}&\hphantom{-}{2}&{-2}&\cr
{-1}&{-2}&\hphantom{-}{2}&\cr
\endpmatrix.
$$
It is recognized that this matrix is the simplest
hyperbolic (i.e., with exactly one negative square)
generalized Cartan matrix of rank $\ge 3$.
It has the smallest possible coefficients for
this type of matrices. There are many publications where the
Kac--Moody algebra $\geg (A_{1,0})$ was considered.
We only mention the  paper of Feingold and Frenkel \cite{FF}.
In that paper there was the first attempt to connect
this algebra with the theory of Siegel modular forms.

The Weyl group of $\geg (A_{1,0})$ is an integral orthogonal
group of signature $(2,1)$, and it is  isomorphic to
the extended modular group $PGL_2(\bz)$. The Weyl-Kac
denominator formula  for $\geg (A_{1,0})$ is
$$
\sum_{g\in PGL_2(\bz)} \hbox{det}(g)
\exp{(2\pi i\,\hbox{tr}((gPg^t-P)Z)}
=
\sum_{\alpha>0}(1-\exp{(2\pi i\,\hbox{tr}(\alpha Z)})
^{{mult}(\alpha)},
\tag{0.1}
$$
where $P=\pmatrix 3&\frac{1}2\\
\frac{1}2&2\endpmatrix$,
$Z=\pmatrix z_1&z_2\\z_2&z_3\endpmatrix$ is a matrix from
the Siegel upper half-plane $\bh$, and the product is taken
over all positive roots of the algebra $\geg (A_{1,0})$
(see \cite{FF}). Unfortunately  an exact formula for
the multiplicities of roots ${\hbox{mult}(\alpha)}$ is not known.
Our point of view on the Lorentzian Kac--Moody algebras theory
(see \cite{N7}, \cite{N8}, \cite{GN1}--\cite{GN3},
where we also use some pioneer ideas due to
R. Borcherds \cite{B1}--\cite{B6}) is that Lorentzian Kac--Moody
algebras of the so called elliptic or parabolic type and
with a lattice Weyl vector (one can consider them as
``the most symmetric") have an automorphic correction.
It means that for a Kac--Moody algebra $\geg (A)$ of
this type, there exists a generalized Kac--Moody
superalgebra without odd real simple roots $\geg (A)^\prime$
such that $\geg (A)\subset \geg (A)^\prime$,
these algebras have the same rank, have the same
real roots system (in particular, any real root is even) and
Weyl group, but $\geg (A)^\prime$ has good automorphic properties
for its Weyl--Kac--Borcherds denominator function $\Phi (z)$.
The last means that $\Phi (z)$ is an automorphic form
on the classical hermitian symmetric domain
of type IV $\Omega$
which one should consider as the complexification
of the hyperbolic space $\La$ where the Weyl group $W$ of $\geg (A)$
acts.
(The Weyl group  of a hyperbolic
Kac--Moody algebra is  a subgroup of an orthogonal
group of signature $(n,1)$.)
The automorphic form  $\Phi (z)$ is invariant with some weight
with respect to an arithmetical orthogonal group
of signature $(n+1,2)$, which contains the Weyl group $W$.

We remark that for a
Lorentzian Kac--Moody algebra $\geg (A)$
its  denominator function  $\Phi (z)$
is invariant only with respect to a hyperbolic
orthogonal group.
(For example, the function \thetag{0.1} does not contain
enough Fourier coefficients to be a Siegel modular form.)
Therefore, one should consider
corrected algebras as
having strictly larger groups of symmetries than
Kac--Moody algebras.
This is the first  advantage
of the notion of automorphic correction.

{}From the point of view of the automorphic form theory,
a form realizing   an automorphic correction is
an automorphic form with a very special Fourier expansion related
with the generalized Cartan matrix $A$
(see \thetag{3.3}, \thetag{3.4}, and \cite{N8}, \cite{GN1} for
the general setting).
Using automorphic properties of
$\Phi (z)$, it is  possible
to find the infinite product expansion of $\Phi (z)$
and calculate multiplicities
of the root spaces decomposition of the corrected algebra
$\geg (A)^\prime$.
This is the second important advantage of these algebras.

Surprisingly,  some well known
classical Siegel modular forms  enable us to construct
automorphic corrections.
In the theory of Siegel modular forms of genus $2$
there is an important  exceptional modular form
$\Delta_5(Z)$of  weight $5$
with non-trivial character $v:\spi\to \{\pm 1\}$.
This function is defined as the product of all even
theta-constants of genus $2$ (see \thetag{1.6} below) and it
delivers
an automorphic correction of the Kac--Moody algebra
$\geg (A_{1,II})$ where
$$
A_{1,II}=
\pmatrix
\hphantom{-}{2}&{-2}&{-2}&\cr
{-2}&\hphantom{-}{2}&{-2}&\cr
{-2}&{-2}&\hphantom{-}{2}&\cr
\endpmatrix
$$
is another simple hyperbolic symmetric
generalized Cartan matrix of rank 3.
This result was proved in \cite{GN1}, \cite{GN2}, where the first
examples of automorphic corrections of elliptic
Lorentzian Kac--Moody algebras were  found.
We  proved   there  that the  following infinite product
expansion is valid
$$
\Delta_5(Z)=
(qrs)^{\frac 1{2}}\,
\prod
\Sb n,\,l,\,m\in \Bbb Z\\
\vspace{0.5\jot}
(n,l,m)>0\,\endSb
\bigl(1-q^n r^l s^m\bigr)^{f(4nm-l^2)}.
\tag{0.2}
$$
The integral exponents $f(4nm-l^2)$ are defined by Fourier
coefficients of a canonical weak Jacobi form of
weight $0$ and index $1$
(see \thetag{1.8}--\thetag{1.9} below).
In \thetag{0.2} we use the variables
$q=\exp(2\pi i z_1)$, $r=\exp(2\pi i z_2)$,
$s=\exp(2\pi i z_3)$.

There exists a   formula for the exponents $f(4nm-l^2)$ in terms of

so-called
H. Cohen numbers. It give us an exact formula for
{\it multiplicities} of root spaces decomposition
of the corrected algebra $\geg (A_{1,II})^\prime$.

Surprisingly,  exactly the same story happens with the
simplest Kac--Moody algebra $\geg (A_{1,0})$.
One of the main results of this paper
is that the automorphic correction of
$\geg (A_{1,0})$ is defined by the first
Siegel modular form $\chi_{35}(Z)$ of odd weight.
This function of   weight $35$
was  constructed by Igusa in \cite{Ig1}
as a linear combination of all ``azygous" triplets of
theta-constants.
Igusa proved that any Siegel modular form
(with the trivial character) of odd weight is divisible
by the form $\chi_{35}$. We call  $\chi_{35}(Z)$
the {\it Igusa modular form}.

To prove that  $\chi_{35}(Z)$ gives an automorphic
correction of the algebra $\geg (A_{1,0})$ and to find the
infinite product  formula for the Igusa
modular form,
we use the representation of  $\chi_{35}(Z)$
in terms of $\Delta_5(Z)=\Delta_5(z_1,z_2,z_3)$.
We show that,
up to a constant, $\chi_{35}(Z)$ is defined by  the quotient
$$
\chi_{35}(Z)=\frac{\chi_{75}(Z)}{\Delta_5(Z)^8}
\tag{0.3}
$$
where  the modular form $\chi_{75}(Z)$ of weight $75$
is defined by the  product
$$
\align
\chi_{75}(Z):=
\kern-10pt\prod_{a,b,c\,mod\, 2}&\kern-8pt
\Delta_5 ({\tsize \frac{z_1+a}2 ,
\tsize\frac{z_2+b}2, \tsize\frac{z_3+c}2})
\prod_{a \,mod\, 2}
\Delta_5 ({\tsize\frac{z_1+a}2},z_2,2z_3)\,
\Delta_5 (2z_1, z_2, {\tsize\frac{z_3+a}2})\\
\vspace{2\jot}
{}&\times\Delta_5 (2z_1, 2z_2, 2z_3)
\prod_{b\,mod\, 2}
\Delta_5 (2z_1, -z_1+z_2, {\tsize\frac{z_1-2z_2+z_3+b}2}).
\tag0.4
\endalign
$$
One can consider the last formula  as
a ``multiplicative" Hecke operator $T(2)$ applied to
$\Delta_5(Z)$. To prove \thetag{0.3}, we compare
divisors of $\chi_{75}(Z)$ and of $\chi_{35}(Z)$ and use the Koecher
principle. At the same time, this gives the new construction of the
form $\chi_{35}(Z)$ and proves that $\chi_{35}(Z)$
is the Siegel modular form of the smallest odd weight.
Due to  \thetag{0.2}--\thetag{0.4}, we  get the infinite product
expansion of $\chi_{35}(Z)$ which is similar to
\thetag{0.2}. Using this expansion, it is
not difficult to prove that $\chi_{35}(Z)$ gives an
automorphic correction of $\geg (A_{1,0})$.

The divisors of the modular forms $\Delta_5(Z)$
and $\chi_{35}(Z)$ are well known.
The divisor of $\Delta_5(Z)$ is
the ``diagonal"
$H_1=\{Z=\pmatrix z_1&0\\0&z_3\endpmatrix\in
\spi\setminus \bh\}$ of the quotient $\spi\setminus \bh$.
The divisor $H_1$, which is the so called
Humbert surface of discriminant $1$, gives moduli of products
of two elliptic curves, and its complement in
$\spi \setminus \bh$ gives moduli of curves of  genus $2$.
(This geometrical description explains applications of
the  modular form  $\Delta_5(Z)$ in the string theory.)
The divisor of the Igusa modular form $\chi_{35}(Z)$
is the sum with multiplicities one of
the Humbert surface $H_1$ above and
the Humbert surface $H_4$ of discriminant $4$.

The consideration with divisors shows us that the quotient
$$
\widetilde{\Delta}_{30}(Z)=\chi_{35}(Z)/\Delta_5(2Z)
$$
is a holomorphic modular form with respect to
the  congruence subgroup
$$
\Gamma_0(2)=\{\pmatrix A&B\\C&D\endpmatrix
\in \spi\ |\ C\equiv 0\mod 2\}.
$$
It turns out that this function  gives an automorphic correction
of the Kac--Moody algebra $\geg (A_{1,I})$
with another very simple generalized Cartan matrix
$$
A_{1,I}=
\pmatrix
\hphantom{-}{2}&{-2}&{-1}&\cr
{-2}&\hphantom{-}{2}&{-1}&\cr
{-1}&{-1}&\hphantom{-}{2}&\cr
\endpmatrix .
$$
Using \thetag{0.2}--\thetag{0.4}, we get the product expansion
of the denominator function
$\widetilde{\Delta}_{30}(Z)$ of $\geg (A_{1,I})^\prime$ and find
multiplicities of its roots.

One can  also consider the second Igusa modular form
$$
\chi_{30}(Z)=\chi_{35}(Z)/\Delta_5(Z)
$$
which is  a Siegel modular form of weight $30$
with the non-trivial character $v: \spi\to \{\pm 1\}$.
The divisor of this function is  $H_4$.
Its infinite product expansion follows from
\thetag{0.2}--\thetag{0.4}. We conjecture that
$\chi_{30}(Z)$ gives an automorphic correction of the
Kac--Moody superalgebra with the generalized Cartan matrix
$A_{1,0}$ and the set of odd indices
$\{2\}\subset \{1,2,3\}$
(see \cite{K2}, \cite{K3} and \cite{R}  for the theory of such
superalgebras).
We don't consider automorphic corrections of Kac--Moody
superalgebras in this paper, but this example shows that they
appear naturally in the subject.
Of course, these considerations show that the four Kac--Moody
algebras which we have considered are very closely related.

Similarly to \thetag{0.4}, applying the Hecke operator
$T(p)$ ($p>2$ is a prime) to $\Delta_5(Z)$, we get
a Siegel modular form $F^{(p^2)}(Z)$ whose  divisor
is the Humbert surface $H_{p^2}$ of discriminant $p^2$, and
we obtain the infinite product and sum expansion of
this form (see Theorem 1.2, Theorem  1.7 and
Appendix A).
We remark that the modular form $F^{(p^2)}(Z)$ is equal
to a finite product (or quotient) of known infinite
series. This leads to satisfactory formulae for
the Fourier coefficients of $F^{(p^2)}(Z)$.

It was proved by van der Geer (see \cite{vdG1}--\cite{vdG2}),
that for any natural $D$,
there exists a Siegel modular form $F^{(D)}(Z)$
whose divisor is the Humbert surface $H_D$.
Thus Theorem 1.7 gives the exact
construction of these modular forms
with divisors $H_{p^2}$.
Moreover, in Appendix B, using the infinite product
formula of the type \thetag{0.2}, we define
Siegel modular forms with divisor $H_D$
for $D$ being not a perfect square.

In Appendix A the general results are obtained
which relate multiplicative Hecke operators acting on
Siegel modular forms,  with Hecke-Jacobi operators
acting on  Jacobi forms.
This is important for the theory of liftings
of Jacobi forms. In particular, we get that
exponents of the product formulae for all modular
forms $F^{(D)}(Z)$ above are related with Fourier
coefficients of Jacobi forms of weight $0$ and index $1$.

In \cite{GN1} and \cite{GN3}
we have pointed out the interpretation of the function
$\Delta_5(Z)$ from the point of view of mirror symmetry for
K3 surfaces and
as the discriminant of K3 surfaces moduli
with the condition $S\subset L_{K3}$ on Picard
lattice of K3 surfaces where $(S)^\perp_{L_{K3}}$ is isomorphic to
the lattice $2U(4)\oplus \langle -2 \rangle$. Here $U$ denote
the even unimodular lattice with signature $(1,1)$.

Similarly, the Igusa form $\chi_{35}(Z)$ serves
for K3 surfaces moduli with  condition
$S=U\oplus E_8(-1)\oplus E_7(-1) \subset L_{K3}$ on
the Picard lattice. Here
$(S)^\perp_{L_{K3}}\cong 2U\oplus \langle -2 \rangle$.

Using automorphic forms $F^{(D)}(Z)$
with the divisor $H_D$,
one can define K3 surfaces submoduli of
codimension one
$\M_{S_1\subset L_{K3}}\subset \M_{S\subset L_{K3}}$, where
$S\subset S_1\subset L_{K3}$ and $\dim S_1 =\dim S +1$,
by the one automorphic form equation.
Here $S\subset L_{K3}$ is one of two conditions on
the Picard lattice of K3 surfaces which we considered above.
This gives K3 surfaces
moduli interpretation of automorphic forms $F^{(D)}(Z)$.
It would be interesting to understand the product and sum
formulae of $F^{(D)}(Z)$ from the point of view of
this algebraic-geometric construction. See \cite{GN3} for
some related considerations.

In Sect. 4 we give the complete list of 13 symmetric hyperbolic
generalized Cartan matrices $A$ of elliptic or parabolic
type and with a lattice Weyl vector such that $A$
contains a submatrix
$
\widetilde{\Bbb A}_1=\pmatrix
\hphantom{-}2 & -2\\
-2&  \hphantom{-}2
\endpmatrix.
$
We hope to publish similar results, as in this paper,
for these $13$ generalized Cartan matrices later. We remark that
the cases $A_{1,II}$ and $A_{2,II}$, from this list,
had been considered in \cite{GN1}, \cite{GN2}, and cases
$A_{1,0}$ and $A_{1,I}$ are considering in this paper.

\head
1. The Hecke operators construction of Siegel modular forms
with rational quadratic divisors and their
infinite product expansions
\endhead

A Siegel modular form of weight $k$ with respect
to $Sp_4(\Bbb Z)$ is,
by definition, a holomorphic function $F(Z)$ on
the Siegel upper half-plane
$$
\Bbb H_2=\{Z=\pmatrix z_1 &z_2\\ z_2 &z_3\endpmatrix
\in M_2(\Bbb C)\ |\ \text{Im\,}(Z)>0\},
$$
that satisfies the functional equation
$$
(F|_k  M)(Z)=F(Z)         \tag{1.1}
$$
where
$$
(F|_k M)(Z):= \hbox{det\,}(CZ+D)^{-k}F
\bigl((AZ+B)(CZ+D)^{-1}\bigr)
$$
for any
$M=\pmatrix A&B\\C&D\endpmatrix\in \spi$.
By $\frak M_k(\spi)$ we denote the finite dimensional space
of all Siegel modular forms of weight $k$.

In this section we define a multiplicative analogue of
Hecke operators which transforms a Siegel modular form
of weight $k$ into another Siegel modular   form
of a different weight.
Let us consider a decomposition  of  a double coset $V$
in the finite union of left cosets with respect to $\spi$
$$
V=\spi M \spi=\sum_i \spi M_i
\tag1.2
$$
where
$M$ and $M_i$ are elements of the group of
symplectic similitudes
$$
GSp_4(\bz)=
\{M\in M_4(\bz)\,|\, {}^tMJM=\mu(M)J,\quad \mu(M)\in \bn \}
\quad ({J}
=\pmatrix 0 & E_2 \\
-E_2 &0\endpmatrix).
$$
\definition{Definition}Let $F(Z)$ be a  Siegel modular form
of weight $k$  and $V$ be the double coset \thetag{1.2}.
We set
$$
[F(Z)]_V:=\prod_i  (F|_k M_i)(Z).
\tag1.3
$$
\enddefinition
The next lemma is evident.
\proclaim{Lemma 1.1} The function $[F(Z)]_V$ on $\bh$
is correctly defined
and  is a Siegel modular form of weight $k\nu$ where $\nu$
is the number of the left cosets  in $V$.
\endproclaim
Let $p$ be a prime.
In accordance with  the  elementary divisors theorem,
there exists only one double coset
$$
T(p)=\spi M\spi=\spi \hbox {diag}(1,1,p,p)\spi
$$
with $\mu(M)=p$.
One can  find a system of representatives
from the distinct  left cosets  in $T(p)$
consisting of $\nu(p)=(p^2+1)(p+1)$ elements:
$$
\align
T&(p)=\spi\pmatrix
p&0&0&0\\0&p&0&0\\0&0&1&0\\0&0&0&1
\endpmatrix+\sum_{a_1,a_2,a_3\,mod\,p}
\spi\pmatrix
1&0&a_1&a_2\\0&1&a_2&a_3\\0&0&p&0\\0&0&0&p
\endpmatrix
\\
\vspace{2\jot}
{}&+
\sum_{a\,mod\,p}
\spi\pmatrix
1&0&a&0\\0&p&0&0\\0&0&p&0\\0&0&0&1
\endpmatrix+
\sum_{b_1,b_2\,mod\,p}
\spi\pmatrix
p&0&0&0\\-b_1&1&0&b_2\\0&0&1&b_1\\0&0&0&p
\endpmatrix.
\tag1.4
\endalign
$$
We denote the modular form $[F(Z)]_{T(p)}$
simply by $[F(Z)]_p$.

\remark{Remark on the definition of $[F(Z)]_p$}
One can give another equivalent definition of the form
$[F(Z)]_p$ using other terms.
Let $M_p=\hbox{diag}(pE_2,E_2)$. Then
$$
\Gamma_0(p)=\spi\cap M_p^{-1}\spi M_p=
\bigl\{\pmatrix A&B\\C&D\endpmatrix\in \spi\,|\,
C\equiv 0\mod p\bigr\}
$$
is a congruence subgroup of the Siegel modular group.
Obviously, one has
$$
T_p=\sum_{\gamma\in \Gamma_0(p)\setminus \spi}
\spi M_p \gamma.
$$
Hence
$$
[F(Z)]_p=\prod_{\gamma \in \Gamma_0(p)\setminus \spi}
F(pZ)|_k\gamma.
$$
The function $F(pZ)$ is a modular form of weight $k$
with respect to $\Gamma_0(p)$, and we may consider
the last product as a multiplicative
$\spi$-symmetrization of a $\Gamma_0(p)$-modular form
(compare with \cite{G2, \S1, \S3} where an additive
symmetrization with respect to a paramodular group
was considered).
\endremark
\medskip

In this section we find $[F(Z)]_p$ for the well-known
Siegel cusp  form  $\Delta_5(Z)$ of weight $5$,
which is defined by the product  of all even theta-constants.
In particular,  we show that the function  $[\Delta_5(Z)]_2$
given  explicitly by the formula  \thetag{0.4}
is connected
with the first non-zero Siegel modular form of odd weight.

We remind first the Igusa's result about the structure of
the graded ring
$$
\frak M_* (Sp_4(\bz))=\bigoplus_{k\in \bn}
\frak M_k (Sp_4(\bz))
$$
of all Siegel modular forms.
This ring has  five generators
$$
\frak M_*  (Sp_4(\bz))=
\bc\,[E_4,\,E_6,\,\chi_{10},\, \chi_{12}, \, \chi_{35}].
$$
The Eisenstein series $E_4(Z)$, $E_6(Z)$ and
the first non-trivial  cusp forms $\chi_{10}(Z)$,
$\chi_{12}(Z)$   (the subscript  denotes  weight)
are algebraically independent. The form
$\chi_{35}(Z)$ is the first Siegel modular form of  odd weight,
which we call the Igusa modular form (see \cite{Ig1}, \cite{Ig2}).

It is known that $\chi_{10}(Z)$ is the square of the  form
$\Delta_5(Z)$ of weight $5$  which is
a $\spi$-modular form of weight $5$ with  non-trivial  character
(see \cite{Fr}).
It means that
$$
(\Delta_5|_5 M)(Z)=
v(M)\,\Delta_5(Z)
\tag1.5
$$
for any $M=\pmatrix A&B\\C&D\endpmatrix\in Sp_4(\bz)$,
where $v:Sp_4(\bz)\to \{\pm 1\}$ is a quadratic character.
Since the commutator subgroup of $\spi$ has index $2$,
there exists only one non-trivial character of $\spi$
(see \cite{M1}, where the exact
formula for $v$ was found for generators of $\spi$).
Up to a constant there exists only one such modular form
of weight $5$.
It can be defined as the product of all even theta-constants
$$
\Delta_5 (Z)
={\tsize\frac 1{64}}\prod_{(a,b)}\vartheta_{a,b}(Z),
\tag1.6
$$
where
$$
\vartheta_{a,b}(Z)=\sum_{l\in \bz^2}
\exp{\bigr(\pi i (Z[l+\frac 1{2}a]+{}^tbl)\bigl)}
\qquad\qquad (Z[l]={}^tlZl)
$$
and the product is taken over all vectors
$a,b\in (\bz/2\bz)^2$
such that ${}^t ab\equiv 0\mod 2$.
(There are exactly ten different $(a,b)$.)
Remark  that $\Delta_5 (Z)$ has integral Fourier coefficients.

We shall apply the operator $[...]_p$ to $\Delta_5(Z)$.
To avoid a problem with the  formal definition,
we fix the system
of representatives  \thetag{1.4} in $T(p)$ taking
$x=0,\dots, p-1$ as a system  modulo  $p$.
Due to \thetag{1.5}  the function $[\Delta_5(Z)]_p$
is a Siegel modular form with  character $v$ or
with  trivial character.
The following result explains our interest to the operator
$[...]_p$.

\proclaim{Theorem 1.2} For $p=2$ the quotient
$$
{[\Delta_5(Z)]_2}/{\Delta_5(Z)^{8}}
$$
is a Siegel modular form of weight $35$ with  trivial character.
For a prime $p>2$ the function
$$
[\Delta_5(Z)]_p/\Delta_5(Z)^{(p+1)^2}
$$
is a Siegel modular form of weight $5p(p^2-1)$ with trivial
character.
\endproclaim
\demo{Proof}We have mentioned that there exists  the only
non-trivial character $v$ of $\spi$ which has  values $\pm 1$.
It is easy to see from \thetag{1.6}   that
$$
v(\nabla)=-1\qquad \text{where }\
\nabla=\pmatrix
1&0&0&0\\0&1&0&1\\0&0&1&0\\0&0&0&1
\endpmatrix.
$$
Let us apply $\nabla$ to $\dtp$. We remind that we fixed
the system of residues $\{0,1,\dots, p-1\}$ modulo $p$.
Using \thetag{1.4}, we obtain that
$$
[\Delta_5(\nabla<Z>)]_p
=v(\nabla)^{p^2+p}v(\nabla^p)^{p+1}\dtp
$$
where $M<Z>=(AZ+B)(CZ+D)^{-1}$ for
$M=\pmatrix A&B\\C&D\endpmatrix$.
The factor $v(\nabla)^{p^2+p}$ equals $1$ for
any  prime $p$.
Thus  $\dtp$ is a Siegel  modular form of weight
$5(p^2+1)(p+1)$ with the trivial character.
To prove the theorem,  we  find the divisor of $\dtp$.

{\it A rational quadratic divisor} of $\bh$ is, by
definition,  the set
$$
\Cal H_\ell=\{\,
\pmatrix z_1&z_2\\ z_2&z_3\endpmatrix\in \bh\ |\
(z_2^2-z_1 z_3)d+cz_3+bz_2+az_1+e=0\,\},
$$
where $\ell=(e, a, b, c, d)\in \bz^5$ is a primitive
(i.e. with the greatest common divisor equals $1$)
integral vector.
The number $D(\ell)=b^2-4de-4ac$ is called
the {\it discriminant} of $\Cal H_\ell$.
This divisor determines the {\it Humbert surface} $H_D$
in the Siegel threefold  $\Cal A_1=\spi\setminus \bh$
(the moduli space of abelian surfaces
with principal polarization)
$$
H_D=\pi\,\bigl(\bigcup
\Sb v\in \bz^5\ primitive\\
\vspace{1\jot}D(v)=D\endSb
\Cal H_v\,\bigr)=\pi\,(\Cal H_\ell)
$$
where $\pi$ is the natural projection $\pi:\bh\to \Cal A_1$.
The following properties of Humbert surfaces in $\Cal A_1$
are well-known and almost evident
if one considers  them from the point  of view of
the orthogonal group $SO(3,2)\cong PSp_4$
(see for example \cite{vdG1}, \cite{GH} where the case of
non-principal polarizations have  also been considered):
\newline
\noindent
${} \qquad \Cal H_\ell$ is not empty $\Leftrightarrow$ $D(l)>0$,
\newline
\noindent
${} \qquad \forall \gamma\in \spi$
the set $\gamma(\Cal H_\ell)$
is a rational quadratic divisor with the same
\newline
\noindent
${} \qquad {}$discriminant $D(\ell)$,
\newline
\noindent
${} \qquad H_D$ is irreducible,
\newline
\noindent
${} \qquad H_D$ can be represented by an equation
$az_1+bz_2+z_3=0$ with $b=0$ or $b=1$.

The next fact  is of the  great importance for different
subjects:
the divisor of $\dtf$ on $\spi\setminus \bh$ is equal to
the Humbert surface
$H_1=\pi(\{Z\in \bh\,|\, z_2=0\})$
of discriminant   $1$.

Let $M_p=\hbox{diag}(p,p,1,1)$, then
$$
M_p(\Cal H_\ell)=\{\,Z\in \bh\ |\
p^2(z_2^2-z_1 z_3)d+pcz_3+pbz_2+paz_1+e=0\,\}.
$$
Thus $D(M_p(\Cal H_\ell))=p^2D(\ell)$ if $(e,p)=1$ and
$D(M_p(\Cal H_\ell))=D(\ell)$ if $(e,p)=p$.
 From the properties mentioned above, it follows
that for any element
$M\in T(p)=\spi \hbox{diag}(p,p,1,1)\spi$
the  divisor
$M(\Cal H_\ell)$
is a rational quadratic divisor of discriminant
$D(\ell)$ or $p^2D(\ell)$.
Therefore the divisor of the modular  form $\dtp$
is a sum of the  Humbert surfaces $H_1$ and  $H_{p^2}$
with some multiplicities $\alpha$ and $\beta$
$$
\hbox{Div}(\dtp)=T(p)^*(H_1)=
\pi\,\bigl(\bigcup
\Sb v\in \bz^5\ primitive\\
\vspace{1\jot}D(v)=1\endSb
\sum_{M_i\in T(p)}M_i^{-1}(\Cal H_v)\,\bigr)
=\alpha H_1+\beta H_{p^2}
$$
where the last sum is taken over representatives $M_i$
from the distinct cosets in $T(p)$.
To define  the integers $\alpha$ (resp. $\beta$),
we have to find the number of the elements $M_i$
from the system $\thetag{1.4}$ such that
$M_i(\Cal H^{(1)})=\Cal H_\ell$
(resp. $M_i(\Cal H^{(p)})=\Cal H_\ell$)
with $D(\ell)=1$ for two quadratic divisors
$$
\Cal H^{(1)}=\{\,Z\in \bh\ |\ z_2=0\,\}
\quad\text{and}\quad
\Cal H^{(p)}=\{\,Z\in \bh\ |\ pz_2=1\,\}.
$$
It is easy to see that $\beta=1$
(there is the only $M_i=\hbox{diag}(p,p,1,1)$)
and $\alpha=(p+1)^2$
(the first and third summands  and the  summands with
$a_2=0$ and $b_1=0$ in \thetag{1.4}).
Hence we proved that
$$
\hbox{Div}(\dtp)=(p+1)^2H_1+H_{p^2}.
$$
Therefore the function
$\dtp/\dtf^{(p+1)^2}$ has no poles on $\bh$.
According to the Koecher principle this function
is holomorphic at infinity.
For $p>2$ the quotient has trivial $\spi$-character.
For $p=2$ we consider $\dtp/\dtf^{(p+1)^2-1}$.
This finishes the proof of the theorem.
\newline
\qed
\enddemo

{}From the  proof of Theorem  1.2, we get

\proclaim{Corollary 1.3}Let $p\ge 2$ be a prime.
The divisor of the Siegel  modular form
\newline
$[\Delta_5(Z)]_p/\Delta_5(Z)^{(p+1)^2}$
of weight $5p(p^2-1)$
coincides with the Humbert surfaces $H_{p^2}$
of discriminant $p^2$ taken with multiplicity one.
\endproclaim

Moreover  we obtain the next results
\proclaim{Corollary 1.4}
1. Any Siegel modular form which has  zero
on $H_{p^2}$ is divisible by the modular form
$\dtp/\dtf^{(p+1)^2}$ of weight $5p(p^2-1)$.

2. (Igusa, \cite{Ig1}) Up to a constant there exists only one
Siegel modular form  $\chi_{35}(Z)$ of weight $35$,
and an arbitrary Siegel modular form
of  odd weigh  is divisible by $\chi_{35}(Z)$.

3. An arbitrary Siegel modular form with the non-trivial character
$v:\spi\to \{\pm 1\}$ can be represented as
$$
\Delta_5(Z)F(Z)\qquad\text{or}\qquad
\frac{[\Delta_5(Z)]_2}{\Delta_5(Z)^9} F(Z)
$$
where $F(Z)$ is a Siegel modular form with the trivial character.
\endproclaim
\demo{Proof} The first statement follows from the Koecher principle
and Corollary 1.3.

Let $F(Z)$ be a Siegel modular form of odd weight.
Then $F(z_1,0,z_3)\equiv 0$ since there do not exist
$SL_2(\bz)$-modular forms of odd weight.
$F(Z)$ is  $0$ on the Humbert surface
$$
H_4=\pi(\{Z=\pmatrix z&z_2\\z_2&z\endpmatrix\in \bh\}),
$$
since $F(z_1,z_2,z_3)=-F(z_3,z_2,z_1)$.
The divisor of the  form
${[\Delta_5(Z)]_2}/{\Delta_5(Z)^8}$
equals $H_1+H_4$.
It proves the second statement.

To prove the third one, we mention that any
$\spi$-modular form
of odd (resp. even) weight with the character  $v$
vanishes on $H_1$ (resp. $H_4$).
\newline
\qed
\enddemo

\remark{Remark}The existence of Siegel modular forms with
the divisors  $H_D$ for any $D$ was proved by van der Geer
using some geometrical considerations
(see \cite{vdG1}, \cite{vdG2}). The forms
$\dtp$ are examples of Siegel modular forms having the divisors
$H_{p^2}$.
In Appendix B  we show how one can construct explicitly
the modular forms with divisors $H_D$ where $D$ is not
a perfect square.
\endremark
\medskip

In the next theorem
we find an infinite product formula for $\dtp$.
To define this formula,  we recall  some  notions and results.

The Fourier-Jacobi expansion of  a Siegel modular form
$F(Z)=F(z_1,z_2,z_3)$
is its Fourier expansion with respect to  $z_3$
$$
F(z_1, z_2, z_3)=f_0(z_1)+
\sum_{m\ge 1}f_{m}(z_1,z_2)\,
\hbox{exp}(2\pi i\,m z_3)
$$
where $z_1=x_1+iy_1\  (y_1>0)$
belongs to the usual upper half-plane  $\Bbb H_1$
and $z_2\in \Bbb C$.
The function
$$
{\tilde f}_{m}(Z):=
f_{m}(z_1,z_2)\,\hbox{exp}(2\pi i\,m z_3)
$$
satisfies  \thetag{1.1} for all elements $M\in \gi$
of the parabolic subgroup
$$
\gi=
\left\{\pmatrix *&0&*&*\\
                *&*&*&*\\
                *&0&*&*\\
                0&0&0&* \endpmatrix
\in Sp_4(\Bbb Z)\right\}.
\tag{1.7}
$$
The Fourier-Jacobi coefficient $f_{m}(z_1,z_2)$
(resp.  ${\tilde f}_{m}(Z)$)
is a Jacobi form (see \cite{EZ} and Appendix A)
of  weight $k$  and index $m$
(resp. a $\gi$-modular form of weight $k$).

The Fourier expansions of the cusp  forms
$\Delta_5(Z)$, $\chi_{10}(Z)$  and  $\chi_{12}(Z)$
can be determined only  by the Fourier coefficients
of their first  Fourier-Jacobi coefficients.
In other words,  they are elements of the so-called Maass subspaces.
For example, the first Fourier-Jacobi
coefficient of $\Delta_5(Z)$ is a Jacobi form of weight $5$
and index $1/2$
$$
\align
\phi_{5,\frac {1}{2}}(z_1,\,z_2)&
=\eta^{9}(z_1)\vartheta_{11}(z_1,z_2)
=-q^{1/2}r^{-1/2}{\kern-2pt}
\prod_{n\ge 1}\,(1-q^{n-1} r)(1-q^n r^{-1})(1-q^n)^{10}\\
{}&=\sum_{n,l\,\equiv 1\,mod\  2} g(n,l)\,
\exp{(\pi i\,(nz_1+lz_2))}.
\endalign
$$
Here  we denote $q=\exp(2\pi i z_1)$ and $r=\exp(2\pi i z_2)$.
Forms $\eta(z_1)$ and $\vartheta_{11}(z_1,z_2)$ are
the Dedekind eta-function and the Jacobi theta-series
respectively.
The Fourier expansion of  $\Delta_5(Z)$
has the  form(see \cite{M2})
$$
\Delta_5(Z)=
\sum\Sb n,l,m\equiv 1\, mod\  2\\
                    n,m>0,\ 4mn-l^2>0\endSb
\sum_{d|(n,l,m)} d^4\, g(\frac{nm}{d^2},\,\frac{l}{d})
\exp{(\pi i\,( nz_1+lz_2+mz_3))}.
$$
Up to a constant  there exists  only one  Jacobi cusp form of
weight $12$ and index $1$. We consider the Jacobi form
of weight $12$ and index $1$ with integral
Fourier coefficients
$$
\align
\phi_{12,1}(z_1, z_2){}&=\frac 1{144}
\bigr(E_4^2(z_1)E_{4,1}(z_1,z_2)-E_6(z_1)E_{6,1}(z_1,z_2)\bigl)\\
\vspace{2\jot}
{}&=(r^{-1}+10+r)q +(10r^{-2}-88r^{-1}-132-88r+10r^{-2})q^2
+\dots\ .
\endalign
$$
Here
$$
E_4(z_1)=1+240 \sum_{n\ge 1}\sigma_3(n)q^n, \qquad
E_6(z_1)=1-504 \sum_{n\ge 1}\sigma_5(n)q^n
$$
are  Eisenstein series for $SL_2(\bz)$,
and $E_{k,1}(z_1,z_2)$
is the  Jacobi--Eisenstein series with integral
Fourier coefficients of  weight $k$  and index $1$
(see \cite{EZ, \S 3}).

Let us introduce a weak Jacobi form
of weight $0$ and index $1$  with integral Fourier coefficients
$$
\align
\phi_{0,1}(z_1, z_2):&=\phi_{12,1}(z_1, z_2)/\Delta_{12}(z_1)=
\sum_{n\ge 0,\, l\in \Bbb Z} f(n,l)\,
\exp{(2\pi i \,(nz_1+lz_2))}\\
{}&=(r^{-1}+10+r)+q(10r^{-2}-64r^{-1}+108-64r+10r^2)+\dots\
\tag{1.8}
\endalign
$$
where
$
\Delta_{12}(z_1)=q\prod_{n\ge 1}(1-q^n)^{24}
$
is the $SL_2(\Bbb Z)$-cusp form of weight $12$.
The Jacobi form  $\phi_{0,1}(z_1, z_2)$
is one of the two canonical generators of the ring
of weak Jacobi forms (see \cite{EZ, \S 9}).
$\phi_{0,1}(z_1, z_2)$ satisfies the same functional
equation as  usual Jacobi forms
(i.e. $\phi_{0,1}(z_1, z_2)\exp{(2\pi i z_3)}$
satisfies \thetag{1.1} for all $M\in \Gamma_\infty$)
and has  nonzero Fourier coefficients
only with indices $(n,l)\in \bz$ such that
$n\ge 0$ (since $\phi_{12,1}(z_1,z_2)$
is a cusp form) and $4n-l^2\ge -1$.
Its weight is even, thus   $f(n,l)=f(n,-l)$ and  $f(n,l)$
depends only on $4n-l^2$.
Therefore we may define a function $f(N)$ of
integral argument such that
$$
f(N)=
\cases
f(n,l)&\  \text{ if } N=4n-l^2\\
0     &\  \text{ otherwise.}
\endcases
\tag1.9
$$
In particular, we have $f(0)=10$, $f(-1)=1$ and $f(n)=0$
if $n<-1$.

In \cite{GN1, Theorem 4.1} we proved that the modular form
$\Delta_5(Z)$ is  the denominator function of
a generalized Lorentzian Kac--Moody
superalgebra and the  product formula
$$
\Delta_5(Z)=
(qrs)^{\frac 1{2}}\,
\prod
\Sb n,\,l,\,m\in \Bbb Z\\
\vspace{0.5\jot}
(n,l,m)>0\,\endSb
\bigl(1-q^n r^l s^m\bigr)^{f(4nm-l^2)}
\tag1.10
$$
is valid
where $q=\exp(2\pi i z_1)$, $r=\exp(2\pi i z_2)$,
$s=\exp(2\pi i z_3)$.
The condition  $(n,l,m)>0$ in the product means
that $n\ge 0$, $m\ge 0$,
$l$ is an  arbitrary integer if $n+m>0$, and
$l<0$ if $n=m=0$.
(We remark that the infinite product absolutely converges for all
$\hbox{Im}(Z)>C$ for a large $C>0$.)

Let us introduce a modular form of weight $35$
$$
\Delta_{35}(Z):=
2^{110}\,\exp(\frac {\pi i}{4})\,
\frac{[\Delta_5(Z)]_2}{\Delta_5(Z)^{8}}.
$$
\proclaim{Theorem 1.5} The Siegel modular form
$\Delta_{35}(Z)$
has integral Fourier coefficients and
can be represented as   infinite product:
$$
\Delta_{35}(Z)=
q^2 r s^2\,(q-s)
\prod\Sb n,\,l,\,m\in \Bbb Z\\
\vspace{0.5\jot}
(n,l,m)>0\,\endSb
\bigl(1-q^n r^l s^m\bigr)^{f_2(4nm-l^2)}.
$$
The integral exponents $f_2(4nm-l^2)$
are defined by the formula
$$
f_2(N)=8f(4N) +2(\biggl(\frac{-N}2\biggl)-1)f(N)
+f(\frac {N}4)
$$
where $f(N)$ is the function \thetag{1.9} and
$$
\biggl(\frac{D}2\biggl)=
\cases
\hphantom{-}1\  \text{ if } D\equiv 1 \mod 8\\
-1\  \text{ if } D\equiv 5 \mod 8\\
\hphantom{-}0\  \text{ if } D\equiv 0 \mod 2.
\endcases
$$
The condition  $(n,l,m)>0$ in the last product means
that $n\ge 0$, $m\ge 0$,
$l$ is an  arbitrary integer if $n+m>0$, and
$l<0$ if $n=m=0$.
\endproclaim

\remark{Remark}The factor $q^2 r s^2\,(q-s)$ shows that
$\Delta_{35}(Z)$ coincides with
$(4i)\chi_{35}(Z)$ where $\chi_{35}$
is the Igusa modular form (see \cite{Ig1, Theorem 3}).
\endremark

Below we give  the proof of Theorem 1.5 together with
its variant for  arbitrary prime $p$.
Before doing this, we formulate  the next corollary.

Let  $\Delta_{30}(Z)=\Delta_{35}(Z)/\Delta_{5}(Z)$
be  a holomorphic cusp form of weight $30$
with non-trivial character $v$
(see the proof of Theorem 1.2).
Using Theorem 1.5 and \thetag{1.10}, we get

\proclaim{Corollary 1.6} The following identity is valid
for $Z$ with a large imaginary part
$$
\Delta_{30}(Z)=
q^{\frac{3}2} r^{\frac{1}2} s^{\frac{3}2}\,(q-s)
\prod\Sb n,\,l,\,m\in \Bbb Z\\
\vspace{0.5\jot}
(n,l,m)>0\,\endSb
\bigl(1-q^n r^l s^m\bigr)^{f'_2(4nm-l^2)}
$$
where the integers $f'_2(4nm-l^2)$
are defined by the formula
$$
f'_2(N)
=8f(4N) +(2\biggl(\frac{-N}2\biggl)-3)f(N)+f(\frac {N}4).
$$
The form $\Delta_{30}(Z)$ has integral Fourier coefficients,
and the divisor of $\Delta_{30}(Z)$ on
$\spi\setminus \bh$ is equal to $H_4$.
\endproclaim
\remark{Remark}Using the product formulae for
$\Delta_{35}(Z)$ and $\Delta_{30}(Z)$, we can calculate the first
non-trivial Fourier-Jacobi coefficients of these functions
having indices $2$ and $3/2$ respectively.
We have
$$
f_2(-4)=f'_2(-4)=1,\quad f_2(-1)=0,\  f'_2(-1)=-1,
\quad f_2(0)=70,\ f_2(0)=60.
$$
Therefore
$$
\phi_{35,2}(z_1,z_2)
=\eta^{69}(z_1)\vartheta_{11}(z_1, 2z_2)
=-q^{2}r^{-1}{\kern-2pt}
\prod_{n\ge 1}\,(1-q^{n-1} r^2)(1-q^n r^{-2})(1-q^n)^{70}
$$
and
$$
\multline
\phi_{30,\frac{3}2}(z_1,z_2)=
\eta^{59}(z_1)\frac{\eta(z_1)\vartheta_{11}(z_1, 2z_2)}
{\vartheta_{11}(z_1, z_2)}=\\
q^{\frac3{2}}r^{-\frac{1}2}
\prod_{n\ge 1}(1+q^{n-1} r)(1-q^{2n-1} r^2)(1-q^{2n-1} r^{-2})
(1+q^n r^{-1})(1-q^n)^{60}.
\endmultline
$$
\endremark
Now we formulate an analogue of Theorem 1.5 for a prime $p>2$.
Let us define a Siegel modular form
$$
F_p(Z):= (-1)^{\frac{p-1}2}p^{5p(2p^2+p+1)}\,
\frac{[\Delta_5(Z)]_p}{\Delta_5(Z)^{(p+1)^2}}
$$
of weight $5p(p^2-1)$.

\proclaim{Theorem 1.7}Let $p$ be  an odd prime.
The form $F_p(Z)$ has the infinite product expansion
$$
F_p(Z)=q^{\frac{5p(p^2-1)}{12}}r^{-\frac{(p-1)}{2}}
s^{\frac{(p^2-1)}{2}}\,
(\frac{r^p-1}{r-1})\,
\prod
\Sb n,\,l,\,m\in \Bbb Z\\
\vspace{0.5\jot}
(n,l,m)\,\endSb
\bigl(1-q^n r^l s^m \bigr)^{f_p(4nm-l^2)}.
$$
The integers $f_p(4nm-l^2)$ are defined by the formula
$$
f_p(N)=p^3f(p^2N) +(p\biggl(\frac{-N}p\biggl)-p-1)f(N)
+f(\frac {N}p)
$$
where
$\biggl(\dsize\frac{\cdot}p\biggl)$
is the Legendre symbol of the quadratic  residue.
The product is taken over  all integral triplets $(n,l,m)$
such that $m\ge 0$, $l$ is an  arbitrary integer
and  $n\ne 0$ or  $m\ne 0$.
In particular,  $F_p(Z)$ has integral Fourier
coefficients.
\endproclaim
\demo{Proof of Theorem 1.5 and Theorem 1.7}Let us consider
the factor
$$
(1-q^n r^l s^m)=(1-\exp{(2\pi i\hbox{tr}(NZ))})
=(1-e(N,Z))
\quad(N=\pmatrix n&l/2\\ l/2&m\endpmatrix),
$$
where $e(N,Z):= \exp{(2\pi i\,\hbox{tr}(NZ))}$.
For any element
$$
M=\pmatrix U&V\\0&p\,{}^tU^{-1}\endpmatrix
\quad\text{with } V=UX,\  X={}^tX\in M_2(\bz),\ \
 UX{\kern-2pt}\mod (p\,{}^tU^{-1}),
$$
from the system of representatives \thetag{1.4}, we have
$$
(1-e(N,\,M<Z>))=(1-e(p^{-1}N[U],\,Z+V)).
$$
Here $N[U]:={}^tUNU$.
Therefore the following formulae are valid:
$$
(1-e(N, Z))|_5 \pmatrix pE_2&0\\0&E_2\endpmatrix
=(1-e(pN, Z))
$$
and
$$
\multline
\prod_{X={}^tX \,mod\,p}(1-e(N, Z))|_5
\pmatrix E_2&X\\0&pE_2\endpmatrix\\
=p^{-10p^3}
\cases
(1-e(N, Z))^{p^2}, &\text{if $N\not\equiv 0 \mod p$}, \\
(1-e(p^{-1}N, Z))^{p^3}, &\text{if $N\equiv 0 \mod p$},
\endcases
\endmultline
$$
for the action of the first and the second summands
from \thetag{1.4}, and
$$
\multline
\prod_{V\, mod\ p{}\,^tU^{-1}}(1-e(N, Z))|_5
\pmatrix U&V\\0&p\,{}^tU^{-1}\endpmatrix\\
= p^{-5p}\cases
(1-e(N[U], Z)), &\text{if $N[U] \not\equiv 0 \mod p$},\\
 (1-e(p^{-1}N[U], Z))^{p}, &\text{if $N[U] \equiv 0 \mod p$},
\endcases
\endmultline
$$
for the action of the third and the fourth summands.

Let us introduce the function $\tilde f: M_2(\bz)\to \bz$
of the matrix argument
$N=\pmatrix n&l/2\\ l/2&m\endpmatrix$
by the formula
$$
\tilde f(N)=\cases f(4\,\hbox{det}\,N) &\text{if }(n,l,m)>0\\
0&\text{otherwise} \endcases
$$
where $f(4\,\hbox{det}\,N)$ is the function \thetag{1.9}.
Collecting  together the above formulae, we obtain
that the factor $(1-e(N, Z))$ appears in the product
$\dtp$ with the exponent
$$
\multline
\tilde f_p(N)=p^3 \tilde f(pN)
+p\sum_{U} \tilde f(p{} N[U^{-1}])+\\
+\cases
\sum_{U} \tilde f( N[U^{-1}])+p^2 \tilde f(N),
&\text{if $N\not\equiv 0 \mod p$}, \\
\tilde f(p^{-1}N), &\text{if $N\equiv 0 \mod p$}.
\endcases
\endmultline
\tag1.11
$$
The last sum is taken over all
$U\in \biggl\{\pmatrix p&0\\0&1\endpmatrix,\,
\pmatrix p&0\\-b&1\endpmatrix\biggr\}_{b=0,\dots, p-1}$
(see \thetag{1.4}).
To find the formula for  $f_p(N)$, one has to calculate
the number (which depends only on $M\mod p$)
of the elements $U$, from the system given above,
such that $M[U^{-1}]$ is an integral matrix.
It provides us  with the given formula for $f_p(N)$.

Due to  \thetag{1.11} and the choice \thetag{1.4},
all factors  $(1-e(N, Z))$  with non-zero $f_p(N)$
have $m\ge 0$.
In the case $p=2$ one easily finds that
the product contains the unique factor with  $n<0$,
namely $(1-q^{-1}s)$.
This  finishs the proof of Theorems 1.5 and 1.7.
\newline
\qed
\enddemo

The formulae of Theorem 1.5 and  Theorem 1.7 show us
that $f_p(N)$ for $p>2$ and $f'_2(N)$ for $p=2$
are the Fourier coefficients of
the   Jacobi form of weight $0$ and index $1$
$$
p^3(\phi_{0,1}|_{0} T^J(p)) (z_1,z_2)-(p+1)\phi_{0,1}(z_1,z_2)
\tag{1.12}
$$
where $T^J(p)$ is the Hecke operator
on the space of Jacobi forms (see \cite{EZ, \S 4}).
We explain this phenomenon in Appendix A
where we give another  ``functorial" proof of
Theorem 1.5 and Theorem 1.7.

\head
2. The orthogonal interpretation of the construction of
$\Delta_{35}$
and application to moduli of K3 surfaces
\endhead

Here we give the simple orthogonal interpretation of
the construction in Sect. 1 of Igusa automorphic forms
using $\Delta_5(Z)$.
Moreover, we give an application to K3 surfaces moduli.

We use notations from \cite{N1} for lattices. For a {\it lattice}
(i.e., an integral symmetric bilinear form)
$S$ and elements $a,b$ of $S$ (i.e., elements of
the module of $S$) we denote by $(a,b)$ the value of
the form of $S$ on elements $a$, $b$.
We denote by $S(q)$ a lattice which one
gets multiplying the form of $S$ by $q\in \bq$.
For an integral symmetric
matrix $A$ we denote by $\langle A \rangle$ a lattice
which is defined by $A$. Thus,
$S=\langle A \rangle$ means that the lattice $S$ has a
bases with the Gram matrix $A$.
We denote by $\oplus$ the orthogonal sum of lattices.

We set
$
U=\left\langle
\matrix
\hphantom{-}0 & -1\\
-1&  \hphantom{-}0
\endmatrix
\right\rangle
$.
The lattice $U$ is unique (up to isomorphism) even unimodular
lattice of signature $(1,1)$.

Let $T$ be a lattice of signature $(t,2)$. One can define
a symmetric domain of type IV $\Omega^+(T)$
which is one of two connected components of
$$
\Omega (T)=\{ Z \in {\Bbb P}(T\otimes \bc) \, \mid \,
(Z, Z)=0,\ (Z, \overline{Z})<0\}.
$$
We also consider the corresponding
homogeneous cones (without zeros)
$\widetilde{\Omega} (T)\subset T\otimes \bc$
and
$\widetilde{\Omega}^+(T) \subset T\otimes \bc$ such that
$\Omega (T)=\widetilde{\Omega} (T)/\bc^\ast$ and
$\Omega^+(T)=\widetilde{\Omega}^+(T)/\bc^\ast$.
Let $O^+(T)$ be the subgroup of index two of $O(T)$ which
keeps the component $\Omega^+ (T)$. It is well-known that
$O^+(T)$ is discrete in $\Omega^+(T)$ and has a fundamental
domain of finite volume.

A function $\Phi$ on $\widetilde{\Omega}^+(T)$ is called
an {\it automorphic form of weight $k$} if $\Phi$ is
holomorphic on $\widetilde{\Omega}^+(T)$,
$\Phi (c\omega)=c^{-k}\Phi (\omega)$ for any $c\in \bc^\ast$,
$\omega \in \widetilde{\Omega}^+(T)$,
and $\Phi (\gamma (\omega))$=
$\chi (\gamma ) \Phi (\omega)$ for any $\gamma \in G$,
$\omega \in \widetilde{\Omega}^+(T)$.
Here $G\subset O(T)^+$ is a subgroup of finite index and
$\chi:G\to \bc^\ast$ some character with the kernel of finite
index in $G$. Then $\Phi$ is called {\it automorphic with
respect to $G$ with the character $\chi$}.
If $\dim \Omega (T)=\rk T-2 \le 2$, the
function $\Phi$ additionally should be holomorphic at infinity of
$\Omega (T)$.
If $\dim \Omega (T)>2$, this condition is
automatically valid by the Koecher principle.

For a lattice $T$ we denote $\Delta^{(2)}(T)=
\{\delta \in T\, \mid \, (\delta , \delta)=2 \}$.
For $e \in T$ with $(e,e) > 0$ we denote
$\Ha_e=\{ Z \in \Omega^+(T)\ \mid \
(Z, e)=0\}$. The $\Ha_e$ is called {\it
the quadratic divisor orthogonal to $e$}. The quadratic
divisor $\Ha_e$ does not change
if one changes $e$ to $te$, $t \in \bq$.

\medskip
{\it The starting point of this paper is
finding in some cases
lattices $T$ of signature $(t,2)$ and
automorphic forms on $\Omega^+(T)$ (with respect to
subgroups of finite index of $O(T)$) with divisor
which is a sum with some multiplicities
of quadratic divisors $\Ha_e$,
$e \in \Delta^{(2)}(T)$.}
For a general discussion see \cite{N9} where similar subject
was considered. Conjecturally, the set of these lattices $T$ is
finite for $\rk T \ge 5$.
\medskip

Let
$$
L_{1,0}=2U\oplus \langle 2 \rangle \ \ \
\text{and}\ \ \
L_{1,II}=2U(4)\oplus \langle 2 \rangle.
$$
Obviously, $L_{1,II} \cong 2L_{1,0}^\ast$ and
$L_{1,0} \cong 2L_{1,II}^\ast$. It
follows that $O(L_{1,0})\cong O(L_{1,II})$.
Using the natural isomorphism
$O^+(L_{1,0})/\{\pm E_5\}\cong Sp_4(\bz)/\{\pm E_4\}$
(see \cite{GN1}, for example), one can identify
Siegel modular  forms (with respect to subgroups of finite
index of $Sp_4(\bz)$)
with automorphic forms on IV type
domains $\Omega^+(L_{1,0})$ or $\Omega^+(L_{1,II})$.

The function $\Delta_5(Z)$ is an automorphic form on
$\Omega^+(L_{1,II})$ with the divisor which is equal
to the sum with multiplicities one of all quadratic
divisors $\Ha_\delta$, $\delta \in \Delta^{(2)}(L_{1,II})$
(see \cite{GN1} where this interpretation of $\Delta_5(Z)$
is given).
By the Koecher principle, this function is unique (up to
multiplying by a constant) and it is automorphic with respect to
$O^+(L_{1,II})$ possibly with some character.
Thus, the function $\Delta_5(Z)$ is defined
(up to a multiplicative constant)
by the lattice $L_{1,II}$ and we denote it by
$\Delta(L_{1,II})_5(Z)$.

Now consider the lattice $L_{1,0}=2U\oplus \langle 2 \rangle$.
We want to construct an automorphic form $\Delta(L_{1,0})$
on $\Omega^+(L_{1,0})$ such that $\Delta(L_{1,0})$ has the divisor
which is equal to the sum with multiplicities one
of all quadratic divisors $\Ha_\delta$,
$\delta \in \Delta^{(2)}(L_{1,0})$.
By the Koecher principle, $\Delta(L_{1,0})$
is unique up to a multiplicative constant if it does exist.

We first remark that elements $\delta \in \Delta^{(2)}(L_{1,0})$
and corresponding quadratic divisors $\Ha_\delta$
are of two different types. An element
$\delta \in \Delta^{(2)}(L_{1,0})$ {\it has type II} if
$(\delta, x)\equiv 0 \mod 2$ for any $x \in L_{1,0}$. Otherwise,
$\delta \in \Delta^{(2)}(L_{1,0})$ {\it has type I};
then there exists $x \in L_{1,0}$ such that
$(x, \delta)\equiv 1 \mod 2$. Two elements of
$\Delta^{(2)}(L_{1,0})$
are conjugate by $O^+(L_{1,0})$ if and only if they have the same
type:
I or II.

Obviously, $\delta \in L_{1,0}$ has type II if and only if
$\delta \in 2L_{1,0}^\ast \cong 2U(4)\oplus \langle 2 \rangle$. Here
$2L_{1,0}^\ast \subset L_{1,0}$ is a canonical sublattice of
$L_{1,0}$
which is equal (or isomorphic) to $2U(4)\oplus \langle 2 \rangle$.
Since the IV type domain of any sublattice of finite index of
$L_{1,0}$ is naturally identified with $\Omega^+(L_{1,0})$,
the automorphic form $\Delta (2L_{1,0}^\ast)_5$ gives the
automorphic form on $\Omega^+(L_{1,0})$
with the divisor which is equal to the
sum with multiplicities one of all quadratic divisors
$\Ha_\delta$, where $\delta \in \Delta^{(2)}(L_{1,0})$ and $\delta$
has type II. Obviously, $O(L_{1,0})=O(2L_{1,0}^\ast)$.
Thus the form $\Delta(2L_{1,0}^\ast)_5$
is automorphic with respect to $O^+(L_{1,0})$ with some character.

Idea of construction of the automorphic form
$\Delta(L_{1,0})$ is using other sublattices
$T\subset L$ such that $T \cong L_{1,II}$,
and considering the product of automorphic forms
$\Delta(T)_5$ of these sublattices.

Let us consider the bases
$(f_1, f_{-1}, f_2, f_{-2}, f_3)$
of $L$ with the Gram matrix
$$
\left(
\matrix
\hphantom{-}0&-1&\hphantom{-}0&\hphantom{-}0&\hphantom{-}0\\
-1&\hphantom{-}0&\hphantom{-}0&\hphantom{-}0&\hphantom{-}0\\
\hphantom{-}0&\hphantom{-}0&\hphantom{-}0&-1&\hphantom{-}0\\
\hphantom{-}0&\hphantom{-}0&-1&\hphantom{-}0&\hphantom{-}0\\
\hphantom{-}0&\hphantom{-}0&\hphantom{-}0&\hphantom{-}0&
\hphantom{-}2\\
\endmatrix
\right).
$$
Obviously, the sublattice
$T_0=\bz 4 f_1+\bz f_{-1}+\bz 2f_2 +\bz 2f_{-2}+\bz f_3 \subset
L_{1,0}$
is isomorphic to $L_{1,II}$. The set $\Delta^{(2)}(T_0)$ contains
some elements of $\Delta^{(2)}(L_{1,0})$ of type I. For example,
$f_{-1}+f_3 \in \Delta^{(2)}(T_0)$ has type I in $L_{1,0}$.
Thus, the automorphic
form $\Delta (T_0)_5$ has the divisor which is equal to
the sum with multiplicities one of all quadratic divisors
$\Ha_\delta$ where
$\delta \in \Delta^{(2)}(L_{1,0})\cap T_0$; and there exists
such $\delta \in \Delta^{(2)}(L_{1,0})\cap T_0$ of type I.

Consider all sublattices $T=g(T_0)\subset L_{1,0}$ where
$g \in O(L_{1,0})$.
Since $4L_{1,0}\subset T_0$, their number is finite. We denote by
$R$ the whole set of these sublattices. Obviously, then
$$
\Psi (L_{1,0})=\prod_{\{T\subset L_{1,0}\}\in R}{\Delta(T)_5}
$$
is an automorphic form on $\Omega^+(L_{1,0})$ with the divisor
which is the sum of quadratic divisors $\Ha_\delta$,
$\delta \in \Delta^{(2)}(L_{1,0})$, with the same multiplicity $a>0$
for all $\delta$ of type I and with the same multiplicity $b \ge 0$
for all $\delta$ of type II. By the Koecher principle, the function
$\Phi (L_{1,0})=\Psi (L_{1,0})\Delta (2L_{1,0}^\ast)_5^{a-b}$ is the
automorphic form on $\Omega^+(L_{1,0})$ with respect to
$O^+(L_{1,0})$ possibly with some character such that
the divisor of $\Phi (L_{1,0})$ is the
sum with the same multiplicity
$a$ for all quadratic divisors $\Ha_\delta$,
$\delta \in \Delta^{(2)}(L_{1,0})$.
It follows that the function
$\Delta (L_{1,0})=\Phi (L_{1,0})^{1/a}$ is the automorphic
form we are looking for.

The orthogonal interpretation of calculations we have done in
Sect. 1 for $p=2$ is that $\# R=15$,
thus, $\Psi (L_{1,0})$ has weight $75$.
The multiplicities $a=1$ and $b=9$.
Thus, $\Delta (L_{1,0})=\Phi (L_{1,0})$ has weight $35$, and further
we denote this automorphic form by $\Delta(L_{1,0})_{35}$.
It is not difficult to repeat these calculations using
the orthogonal language. Discriminant forms
technique (see \cite{N1}) is very useful for these
calculations. We leave these calculations to an
interesting reader.

Now it is not difficult to identify $\Delta (L_{1,0})_{35}$ with
the Igusa modular  form $\chi_{35}$ since these forms
have the same odd weight and $\chi_{35}$ is automorphic
with respect to $O^+(L_{1,0})$, i.e., it
has the trivial character on $SO^+(L_{1,0})$; it follows it has
the character $\det (g)$ for $g \in O^+(L_{1,0})$ since
the weight is odd.

Consider $\delta \in \Delta^{(2)}(L_{1,0})$. This
$\delta$ defines the reflection $s_\delta \in O^+(L_{1,0})$
where
$$
s_\delta(x)=x-(\delta,x)\delta,\ \  x \in L_{1,0}.
$$
It follows that $s_\delta(\delta)=-\delta$ and $s_\delta (x)=x$
if $(\delta, x)=0$. Thus, $\det (s_\delta)=-1$.
We have $\chi_{35}(\omega)=
\chi_{35}(s_\delta(\omega))=-\chi_{35}(\omega)=0$
if $\bc \omega \in \Ha_\delta $.
Thus $\chi_{35}$ is equal to zero on
the quadratic divisor $\Ha_\delta$.
By the Koecher principle,
$\chi_{35}/\Delta (L_{1,0})_{35}$ is a
holomorphic automorphic form of weight $0$. Thus, it is
a constant and $\Delta (L_{1,0})_{35} = c\chi_{35}$
where $c \in \bc^\ast$.

We mention the following application of these considerations to
geometry of moduli of K3 surfaces
(see \cite{N9} for a general setting).

\proclaim{Theorem 2.1} Igusa modular form
$\chi_{35}=c\Delta(L_{1,0})_{35}$, $c \in \bc^\ast$,
gives the discriminant of K3 surfaces moduli
with condition
$S=U\oplus E_8(-1)\oplus E_7(-1)\subset L_{K3}$
on Picard lattice where $E_8$ and $E_7$ are the  standard
positive definite lattices which are defined
by the  root systems of the same type.
\endproclaim
\demo{Proof} The lattice $L_{K3}=3U \oplus 2E_8(-1)$ is the
even unimodular lattice of signature $(3,19)$.
Using discriminant form technique (see \cite{N1}),
it is very easy to calculate  that
$T=(S)^\perp_{L_{K3}} \cong 2U \oplus \langle -2
\rangle=L_{1,0}(-1)$.
The discriminant $\widetilde{D}(S\subset L_{K3})$
of K3 surfaces with condition
$S\subset L_{K3}$ on Picard lattice is equal to union
of all quadratic divisors $\Ha_{\delta}$, $\delta \in
\Delta^{(-2)}(T)=\Delta^{(2)}(L_{1,0})$.
Thus, Igusa modular form $\chi_{35}$ is equal to zero
exactly on this discriminant. It follows the statement.
\newline
\qed
\enddemo

Similarly, $\Delta_5(Z)$ gives the discriminant of K3 surfaces moduli

with condition $S\subset L_{K3}$ on Picard lattice where
$S$ is the hyperbolic lattice with signature $(1,16)$ and
with the discriminant quadratic form
$2u_+^{(2)}(4)\oplus q_{1}^{(2)}(2)$. For this condition,
$(S)^\perp_{L_{K3}}\cong 2U(4)\oplus \langle -2 \rangle =
L_{II}(-1)$. See \cite{GN1} and
\cite{GN3} where this case was considered.

At the same time,
$\Delta_5(Z)=\Delta (2L_{1,0}^\ast)_5(Z)$
is equal to zero on the type II
part of the discriminant of K3 surfaces moduli
with condition $U\oplus E_8(-1)\oplus E_7(-1)\subset L_{K3}$
on Picard lattice,  which is union of all quadratic divisors
$\Ha_\delta$, $\delta \in \Delta^{(2)}(L_{1,0})$
and $\delta$ has type II.

The Igusa modular  form $\chi_{30}=c_1\chi_{35}/\Delta_5=
c_2\Delta (L_{1,0})_{35}/\Delta (2L_{1,0}^\ast)_5$,
$c_1,c_2 \in \bc^\ast$,
on $\Omega^+(L_{1,0})$ gives the type I part of this discriminant
since its divisor is the sum (with
multiplicities one) of all quadratic divisors $\Ha_\delta$,
where $\delta \in \Delta^{(2)}(L_{1,0})$ and $\delta$ has type I.
Let $L_{1,I}$ be a sublattice of $L_{1,0}$ which is generated by all
elements $\delta \in \Delta^{(2)}(L_{1,0})$ of
type I. One can easily find that $L_{1,I}=L_{1,0}$. It follows that
$\chi_{30}$ is not equal to the discriminant (or the full
discriminant if one likes) of
K3 surfaces moduli with condition on Picard lattice.

\medskip

Let us consider a condition
$S=U\oplus E_8(-1)\oplus E_7(-1)\subset L_{K3}$
on Picard lattice of K3 which we have considered above.
Any primitive intermediate hyperbolic sublattice
$S\subset S_1 \subset L_{K3}$ defines a condition $S_1\subset
L_{K3}$
on Picard lattice of K3 surfaces. It defines the submoduli
$\M_{S_1\subset L_{K3}}\subset  \M_{S\subset L_{K3}}$ of
moduli $\M_{S\subset L_{K3}}$ of K3 surfaces with the
condition $S\subset L_{K3}$ on Picard lattice.
Consider a negative definite sublattice $Q=S_1\cap L \subset L$
where $L=(S)^\perp_{L_{K3}}$.
If $Q=\bz \delta$ has the rank one and $\delta^2=-2p^2$,
the submoduli $\M_{S_1\subset L_{K3}}$ are the set of zeros
of the Siegel modular form $F_p(Z)$ which we have
constructed in Sect. 1.
Similarly, one can use modular
forms constructed in Sect. 1 and Appendix B
to define any codimension one submoduli
$\M_{S_1\subset L_{K3}}\subset \M_{S\subset L_{K3}}$
by one modular form (with respect to
the full group $O(L)^+$) equation.
It gives the K3 surfaces interpretation of  these
Siegel modular forms.
It is interesting to understand the algebraic-geometric
sense of the infinite product (and the sum as well) expansions
which we have found in Sect. 1 and Appendix B,
of the corresponding modular forms.
Compare with some considerations in \cite{GN3}.

\vskip5pt

Similar results we have for the condition
$S\subset L_{K3}$ when
$S^\perp_{L_{K3}}=2U(4)\oplus \langle -2 \rangle$.

\vskip5pt

For the automorphic correction of Lorentzian Kac--Moody algebras
which
we consider below, it is important to have automorphic forms
$\Phi$ on $\Omega (T)$ with the divisor which is a sum of some
quadratic divisors $\Ha_\delta$,
$\delta \in \Delta^{(2)}(T)$, with the multiplicity one.
The automorphic forms $\Delta(L_{1,0})_{35}$ and
$\Delta(L_{1,II})_5$
give such automorphic forms. Quotients
$\Phi (L_{1,0})_{30}=\Delta(L_{1,0})_{35}/\Delta (T)_5$,
where $T\subset L_{1,0}$ and $T\cong L_{1,II}$ also give
automorphic forms with these property by the Koecher principle.
There are three different types of embeddings
$T\subset L_{1,0}$ with $T\cong L_{1,II}$.

The first type is defined
by the condition
$L_{1,0}/T \cong (\bz/2\bz)^4$. It gives $T=2L_{1,0}^\ast$ and
$\Delta_{30}=\Delta(L_{1,0})_{30}=
\Delta(L_{1,0})_{35}/\Delta(2L_{1,0}^\ast)_5=c\chi_{30}$.

The second type  $\{T\subset L_{1,II}\}\in R$ we have considered
above. For this type, $L_{1,II}/T\cong (\bz/4\bz)\dotplus
(\bz/2\bz)^2$.
It gives $\widetilde{\Delta}_{30}=
\widetilde{\Delta}(L_{1,0})_{30}=\Delta(L_{1,0})_{35}/\Delta(T)_5$.

There is the third type when
$T = \bz f_1 + \bz 4f_{-1}+ \bz f_2 + \bz 4f_{-2}+\bz f_3
\subset L_{1,0}$ for the bases of $L_{1,0}$ above. For
this type $L_{1,0}/T\cong (\bz/4\bz)^2$.
It gives
$\widetilde{\widetilde{\Delta}}_{30}=
\widetilde{\widetilde{\Delta}}(L_{1,0})_{30}=
\Delta(L_{1,0})_{35}/\Delta(T)_5$.

In  Sect. 3 we use Fourier expansions at the appropriate cusps of
the IV type domain $\Omega^+(L_{1,0})$ (a cusp is
defined by the  primitive isotropic element of $L_{1,0}$)
of the automorphic forms
$\Delta_5=\Delta(L_{1,II})_5$, $\Delta_{35}=\Delta (L_{1,0})_{35}$
and $\widetilde{\Delta}_{30}=\widetilde\Delta(L_{1,0})_{30}$
for the automorphic correction
of some Lorentzian Kac--Moody Lie algebras.
It seems that the Igusa form
$\Delta_{30}=\Delta(L_{1,0})_{30}=c\chi_{30}$ and the automorphic
form
$\widetilde{\widetilde{\Delta}}_{30}=
\widetilde{\widetilde{\Delta}}(L_{1,0})_{30}$ may be used for the
automorphic correction of some Lorentzian Kac--Moody Lie
superalgebras
(see Remark 3.2 below), but we don't consider automorphic correction
for this more general class of Lorentzian Kac--Moody algebras in
this paper.

\head
3. Automorphic correction of ``the simplest''
Lorentzian Kac--Moody algebras or rank $3$
using Igusa automorphic forms
\endhead

There are three closely related elliptic
hyperbolic (here we use terminology of \cite{N8})
generalized
Cartan matrices of the rank $3$ (see \cite{K1}
for general definitions and results related
with Kac--Moody algebras):
$$
A_{1,0}=
\pmatrix
\hphantom{-}{2}&\hphantom{-}{0}&{-1}&\cr
\hphantom{-}{0}&\hphantom{-}{2}&{-2}&\cr
{-1}&{-2}&\hphantom{-}{2}&\cr
\endpmatrix, \qquad
A_{1,I}=
\pmatrix
\hphantom{-}{2}&{-2}&{-1}&\cr
{-2}&\hphantom{-}{2}&{-1}&\cr
{-1}&{-1}&\hphantom{-}{2}&\cr
\endpmatrix ,
$$
$$
A_{1,II}=
\pmatrix
\hphantom{-}{2}&{-2}&{-2}&\cr
{-2}&\hphantom{-}{2}&{-2}&\cr
{-2}&{-2}&\hphantom{-}{2}&\cr
\endpmatrix.\qquad{}
$$
They (especially the first one $A_{1,0}$) give the simplest
generalized Cartan matrices of this type. It does not mean that
the corresponding to these matrices Lorentzian
Kac--Moody algebras are the simplest ones.
It is why we use `` \ '' in the title of the
paper and this Section.

We consider a hyperbolic lattice $M_{1,0}=\bz f_2 + \bz f_{-2}
+\bz f_3$ with the symmetric bilinear form
$(f_2,f_2)=(f_{-2},f_{-2})=0$,
$(f_2,f_{-2})=-1$, $(f_3,f_3)=2$, $(f_2,f_3)=(f_{-2},f_3)=0$.
Thus, the Gram matrix of $f_2,f_{-2},f_3$ is equal to
$$
\left(
\matrix
\hphantom{-}0 &-1 & \hphantom{-}0 \\
-1& \hphantom{-}0 & \hphantom{-}0\\
\hphantom{-}0 &\hphantom{-} 0 & \hphantom{-}2
\endmatrix
\right).
$$
First, we show that the generalized Cartan matrices
$A_{1,0}$, $A_{1,I}$ and $A_{1,II}$
are defined by the lattice $M_{1,0}$ and some its
natural sublattices which we introduce below (here we repeat
some results of \cite{N5} and \cite{GN1}).
The lattice $M_{1,0}$ is
hyperbolic, i.e. it has signature $(2,1)$. Thus $M_{1,0}$
defines a cone
$$
V(M_{1,0})=\{x \in M_{1,0}\otimes \br\ | \ (x,x)<0\}
$$
which is union of its two half-cones.
(All general definitions and notations below hold
for an arbitrary hyperbolic (i.e. with signature $(n,1)$)
lattice, and we shall use them later.)
We choose one of this half-cones $V^+(M_{1,0})$.
We denote by $O^+(M_{1,0})$ the subgroup of $O(M_{1,0})$
of index two which fixes
the half-cone $V^+(M_{1,0})$. It is well-known that
the group $O^+(M_{1,0})$ is discrete in the corresponding hyperbolic
space $\La^+(M_{1,0})=V^+(M_{1,0})/\br_{++}$ and has
a fundamental domain of finite
volume (below indices $++$ and $+$ denote positive and
non-negative numbers respectively).
Any reflection $s_\delta \in O(M_{1,0})$ with respect to an
element $\delta \in M_{1,0}$ with $(\delta, \delta)>0$ is a
reflection in the hyperplane
$$
{\Cal H}_\delta =\{ \br_{++}x \in \La^+(M_{1,0})\ |\ (x, \delta)=0\}
$$
of $\La^+(M_{1,0})$. This maps the half-space
$$
 \Ha^+_\delta =\{ \br_{++}x \in \La^+(M_{1,0})\ |
\ (x, \delta)\le 0 \}
$$
to the opposite half-space $\Ha^+_{-\delta}$ which are both bounded

by the  hyperplane ${\Cal H}_\delta$. Here $\delta \in M_{1,0}$ is
called orthogonal to $\Ha_\delta$ and $\Ha^+_\delta$.
All reflections of $M_{1,0}$ generate the
reflection subgroup $W(M_{1,0})\subset O^+(M_{1,0})$.

\smallpagebreak

The hyperbolic lattice $M_{1,0}$
is very special, and its automorphism group is
well-known (see \cite{N5}, for example).
Firstly, $O^+(M_{1,0})=W^{(2)}(M_{1,0})$
where  index $2$ denote the subgroup
generated by reflections in all elements of $M_{1,0}$ with square
$2$.
Thus, $O^+(M_{1,0})$ is generated by
reflections in $\Delta^{(2)}(M_{1,0})$.
An element $\delta \in \Delta^{(2)}(M_{1,0})$ and the
corresponding reflection $s_\delta$ have one of two types:

\vskip5pt

Type ${}\hphantom{I}$I: $(\delta, M_{1,0})= ${}\hphantom{2}$\bz$.

Type II: $(\delta, M_{1,0})=2\bz$.

\vskip5pt

We introduce sublattices $M_{1,I}$ and $M_{1,II}$ which are generated

by
all elements $\delta \in \Delta^{(2)} (M_{1,0})$ of the type $I$ and
$II$
respectively.
We have
$$
\align
M_{1,I}&
=\{ nf_2+lf_3+mf_{-2} \in M_{1,0} \ | \ n+l+m \equiv 0\mod 2\},\\
M_{1,II}&
=\{nf_2+lf_3+mf_{-2} \in M_{1,0}\ |\ n \equiv m \equiv 0\mod 2\}.
\endalign
$$
The lattice $M_{1,I}$ is unique in its genus and is defined
by its discriminant form $q_{M_{1,I}}\cong q^{(2)}_{5}(8)$
(we use notations from \cite{N1}). There exists a
unique overlattice $M_{1,I}\subset M$ of index $2$, and it is equal
to $M_{1,0}$. The lattice
$M_{1,II}=2M_{1,0}^\ast \cong U(4)\oplus <2>$. The lattice
$M_{1,0}=2M_{1,II}^\ast $. It follows that sublattices or
overlattices $M_{1,0}$, $M_{1,I}$, $M_{1,II}$
are defined by one another, and their automorphism groups
are naturally identified.

An element $\delta\in \Delta^{(2)} (M_{1,0})$ has the type I
(resp. II) if and only if
$\delta \in M_{1,I}$ (resp. $\delta \in M_{1,II}$). It follows that
the subgroup of $O^+(M_{1,0})=W^{(2)}(M_{1,0})$
generated by all reflections
$s_\delta$ of the type I (resp. II) is equal to
$W^{(2)}(M_{1,I})$
(resp. $W^{(2)}(M_{1,II})$).
Obviously, three sublattices $M_{1,0}$, $M_{1,I}$, $M_{1,II}$ are
$W^{(2)}(M_{1,0})$-invariant, and
both subgroups  $W^{(2)}(M_{1,I})$ and $W^{(2)}(M_{1,II})$ are
normal
in $W^{(2)}(M_{1,0})$.
The index $[W^{(2)}(M_{1,0}):W^{(2)}(M_{1,I})]=2$ and
$[W^{(2)}(M_{1,0}):W^{(2)}(M_{1,II})]=6$. Fundamental polyhedra
$\M_{1,0}$, $\M_{1,I}$, and $\M_{1,II}$ for groups
$W^{(2)}(M_{1,0})$, $W^{(2)}(M_{1,I})$ and $W^{(2)}(M_{1,II})$
respectively are equal to
$\M_{1,i}=\bigcap_{\delta \in P(\M_{1,i})}\Ha_\delta^+$
where $P(\M_{1,i})$ are
minimal (with this equality) sets of elements of
$\Delta^{(2)}(M_{1,i})$, where $i=0,I, II$. They are
called orthogonal vectors to polyhedra $\M_{1,i}$ and are
equal to
$$\align
P(\M_{1,0})&=\{-f_2+f_{-2},\ f_3,\ f_2-f_3\};\\
P(\M_{1,I})&=\{f_2+f_3,\ f_2-f_3,\ -f_2+f_{-2}\};\\
P(\M_{1,II})&=\{2f_2-f_3,\  2f_{-2}-f_3,\  f_3\}.
\endalign
$$
Here $\M_{1,0}$ is a triangle with angles $\pi/2$, $0$, $\pi/3$;
$\M_{1,I}$ is a triangle with angles $\pi/3$, $0$, $\pi/3$ and
$\M_{1,II}$ is a triangle with zero angles (i.e., with its
vertices at infinity).
It follows that the groups $W^{(2)}(M_{1,0})$, $W^{(2)}(M_{1,I})$,
$W^{(2)}(M_{1,II})$ are generated by reflections
in elements of $P(\M_{1,0})$, $P(\M_{1,I})$ and
$P(\M_{1,II})$ respectively. We denote by
$$
\text{Sym\ }(P(\M_{1,i}))=
\{g \in O^+(M_{1,i})\ |\ g(P(\M_{1,i}))=P(\M_{1,i}) \}
$$
the {\it group of symmetries} of the fundamental
polyhedron $\M_{1,i}$
and its set $P(\M_{1,i})$ of orthogonal vectors.
The group $\text{Sym\ }(P(\M_{1,0}))$ is trivial,
$\text{Sym\ }(P(\M_{1,I}))$
has order two and is
generated by $s_{f_3}$, the group
$\text{Sym\ }(P(\M_{1,II}))$ is the
group of symmetries of the right triangle (i.e. it is $S_3$) and
is generated by $s_{f_2-f_3}, s_{f_{-2}-f_2}$.
We can write down the groups
$O^+(M_{1,0})=O^+(M_{1,I})=O^+(M_{1,II})$
as the semi-direct products:
$$
\split
O^+(M_{1,0})&=O^+(M_{1,I})=O^+(M_{1,II})=W^{(2)}(M_{1,0})\\
&=W^{(2)}(M_{1,I})\rtimes\text{Sym\ }(P(\M_{1,I}))
=W^{(2)}(M_{1,II})\rtimes \text{Sym\ }(P(\M_{1,II})).
\endsplit
$$
The Gram matrices of sets $P(\M_{1,0})$,
$P(\M_{1,I})$ and $P(\M_{1,II})$ are respectively
equal to the generalized Cartan matrices
$A_{1,0}$, $A_{1,I}$ and $A_{1,II}$ above. It shows that
these matrices are naturally related with the lattice $M_{1,0}$
and its natural sublattices $M_{1,I}$ and $M_{1,II}$.

To unify notation, further $i=0$, $I$ or $II$. The lattice
$M$ denote $M_{1,i}$, and $\M$ denote
the fundamental polygon $\M_{1,i}$ above for $W=W^{(2)}(M)$,
and $P(\M)=\{\delta_1,\delta_2,\delta_3\}$.
The generalized Cartan matrices $A=A_{1,i}$(equivalently the Weyl
groups
$W^{(2)}(M)$), have so called {\it elliptic type}. This means
that the fundamental polygon $\M$ has finite volume in $\La^+(M)$.
Since the cone $V^+(M)$ is self-dual, it
is equivalent to any of the embeddings
$$
\br_+\M \subset \overline{V^+(M)}\subset
\br_+\delta_1+\br_+\delta_2 +\br_+\delta_3
\tag{3.1}
$$
where
$$
\br_+\M=(\br_+\delta_1+\br_+\delta_2 +\br_+\delta_3)^\ast =
\{x \in M\otimes \br\ |\  (x, P(\M))\le 0\}.
$$
Another very important property of $A_{1,i}$ is that
$P(\M_{1,i})$ has a {\it lattice Weyl vector} $\rho=\rho_i$ which
is,
by definition, an element
$\rho_i \in M_{1,i}\otimes \bq$ such that
$$
(\rho_i, \delta )=(\delta, \delta )/2=-1\
\text{for any}\ \delta \in P(\M_{1,i}).
$$
We have
$$
\rho_0=3f_2+2f_{-2}-f_3/2,\ \ \rho_I=2f_2+f_1,\ \
\rho_{II}=f_2+f_{-2}-f_3/2.
\tag3.2
$$

Consider the {\it complexified cone}
$\Omega (V^+(M))=M\otimes \br+iV^+(M)$. Remark that
for all our lattices $M=M_{1,i}$ these cones are identified
(by the isomorphism of lattices over $\bq$), and
we always can use the same coordinates $(z_1,z_2,z_3)$ for the
point $z=z_3f_2+z_2f_3+z_1f_{-2}$ of these domains.
Considering the lattice $L_k=U(k)\oplus M$ where
$U(k)=
\left(\matrix
0 &-k\\
-k &0
\endmatrix
\right)$, $k\in \bn$, we identify
$\Omega(V^+(M))$ with one (of two) connected component of
the symmetric domain
$$
\Omega(L_k)=\{Z \in {\Bbb P}(L_k\otimes \bc)\ |\ (Z,Z)=0,\
(Z,\overline{Z})<0\}
$$
of type IV. If $e_1,e_{-1}$ is the bases of $U(k)$,
one should correspond to
$z \in \Omega (M)$ the point
$\bc (((z,z)/2)e_1+(1/k)e_2+z)\in \Omega (L_k)$.Using these
embeddings, we
can speak about {\it automorphic forms} on
$\Omega (V^+(M))$ with respect to subgroups of finite index of
$O(L_k)$. Thus, an {\it automorphic form} on $\Omega (M)$ means
a holomorphic function on $\Omega (V^+(M))$ which is an automorphic
form (with some non-negative weight) with respect to some subgroup
of finite index of $O(L_k)$ for some $k$.

One can correspond to the matrix $A=A_{1,i}$ the Kac--Moody algebra
$\geg=\geg (A)$ (see \cite{K1}).
{\it An automorphic correction} (see \cite{B3}, \cite{B4}
and \cite{N8} and \cite{GN1} for a general setting)
of the Kac--Moody algebra
$\geg (A)$ is an automorphic form $\Phi (z)$ on
$\Omega (V^+(M))$, $M=M_{1,i}$ which has the Fourier
expansion
$$
\Phi (z)=
\hskip-2pt
\sum_{w\in W}{ \det(w)
\left(\exp(-2\pi i (w(\rho),z))
- \hskip-14pt  \sum_{ a \in M^\ast \cap \br_{++}\M}
\hskip-4pt{m(a)\exp(-2\pi i (w(\rho+a),z))}\right)},
\tag3.3
$$
where all $m(a)$ are integers.
(This is  a very special and very restricted type of
Fourier expansion.)
Using this Fourier expansion of $\Phi (z)$, one can
construct a {\it generalized Kac--Moody Lie superalgebra without odd
real simple roots}
$\geg(M, {}_s\Delta)$ which has the Weyl--Kac--Borcherds denominator
function $\Phi (z)$. This algebra contains the Kac--Moody algebra
$\geg(A)$ and has better automorphic properties for its denominator
function. Using these automorphic properties it
is possible to calculate the product formula for $\Phi (z)$
(the product converges for $z$ with large $\text{Im}~z \in V^+(M)$):
$$
\Phi (z)=
\sum_{w\in W}{ \det(w)
\left(\exp(-2\pi i (w(\rho),z))
- \hskip-18pt  \sum_{ a \in M^\ast \cap \br_{++}\M}
{m(a)\exp(-2\pi i (w(\rho+a),z))}\right)}
$$
$$
=\exp{\left(-2\pi i(\rho,z)\right)}
\prod_{\alpha \in \Delta(M)^{\re}_+\cup \Delta(M^\ast)^{\im}_+}
{\left( 1-\exp{ \left(-2\pi i (\alpha ,
z)\right)}\right)^{\mult~\alpha}}.
\tag3.4
$$
Here $\Delta(M)^{\re}_+=(W(P(\M))\cap Q_+)$ where
$Q_+=\bz_+ \delta_1+\bz_+ \delta_2 + \bz_+ \delta_3$;
equivalently, $\Delta(M)^{\re}_+=\{\delta \in \Delta^{(2)}(M)|
(\rho, \delta)<0\}$. For an intermediate
sublattice $M\subset T \subset M^\ast$,  we denote
$\Delta(T)^{\im}_+=\overline{V^+(M)}\cap T$.
Elements of
$\Delta(M)^{\re}_+$ are called {\it positive real roots}
(they have the square $2$).
Elements of $\Delta(T)^{\im}_+$ are called
{\it positive imaginary roots} (they have non-positive squares).
The exponents $\mult~\alpha$ are integers. They are called
{\it multiplicities} of roots $\alpha$ and are very
important invariants of  the Kac--Moody superalgebra $\geg (M,
{}_s\Delta)$.
Always $\mult~\alpha =1$
if $\alpha \in \Delta(M)^{\re}_+$. We refer to \cite{GN1} for
details.

Using results of Sect. 1, we find automorphic corrections
for Kac--Moody algebras $\geg (A)$, $A=A_{1,0}$, $A_{1,I}$,
$A_{1,II}$. For the algebra $A_{1,II}$ this have been done in
\cite{GN1}, \cite{GN2} using $\Delta_5(Z)$.

We identify a point
$Z=\left(\matrix z_1&z_2\\ z_2&z_3\endmatrix\right)\in
{\Bbb H}_2$ with the point
$z=z_3f_2+z_2f_3+z_1f_{-2} \in \Omega (V^+(M))$
with the coordinates
$(z_1,z_2,z_3)$ which we have introduced above.
Then Siegel automorphic forms on ${\Bbb H}_2$
will give some automorphic
forms on $\Omega (V^+(M))$ because of
the well-known isomorphism
$Sp_4(\bz)/\{\pm E_4\}\cong O^+(L_1)/\{\pm E_5\}$ (see
\cite{GN1}, for example).
Thus, we can consider the automorphic form
$\Delta_5(Z)$ and the
Igusa automorphic forms
$\chi_{35}(Z)$ and $\chi_{30}(Z)$ which we have studied
in Sect. 1, as automorphic forms on $\Omega (V^+(M))$.
We can rewrite
$$
\exp{(\pi i(nz_1+lz_2+mz_3))}=
\exp{(-\pi i (a,z))},\ \ a=nf_2-lf_3/2+mf_{-2}
$$
where $z=z_3f_2+z_2f_3+z_1f_{-2}$.

Using the lattice $M_{1,II}$, we can
rewrite $\Delta_5(Z)$ as
$$
\Delta_5(z)
=\exp{\left(-\pi i(\rho_{II},z)\right)}
\prod_{\alpha \in \Delta(M_{1,II})^{\re}_+\cup
\Delta(M_{1,II})^{\im}_+}
{\left( 1-\exp{ \left(-\pi i (\alpha , z)\right)}\right)^
{\mult~\alpha}},
\tag3.5
$$
where
$$
\mult~\alpha = f(-(\alpha,\alpha)/2)
\tag3.6
$$
for the function $f(N)$ introduced in \thetag{1.9}.
It follows that
$$
\Delta_5(2z)
=\exp{\left(-2\pi i(\rho_{II},z)\right)}
\prod_{\alpha \in \Delta(M_{1,II})^{\re}_+\cup
\Delta(M_{1,II})^{\im}_+}
{\left( 1-\exp{ \left(-2\pi i (\alpha , z)\right)}\right)^
{\mult~\alpha}},
\tag3.7
$$
$\mult~\alpha = f(-(\alpha, \alpha)/2)$, has the form \thetag{3.4}.

Using the lattice $M_{1,0}$, we can rewrite $\Delta_{35}(z)$ as
$$
\Delta_{35}(z)
=\exp{\left(-2\pi i(\rho_{0},z)\right)}
\prod_{\alpha \in \Delta(M_{1,0})^{\re}_+\cup \Delta(M_{1,0}^\ast
)^{\im}_+}
{\left( 1-\exp{ \left(-2\pi i (\alpha ,
z)\right)}\right)^{\mult~\alpha}},
\tag3.8
$$
where
$$
\mult~\alpha = f_2(-2(\alpha, \alpha))
\tag3.9
$$
for the function $f_2(N)$ introduced in Theorem 1.5.

We consider an automorphic form
$\widetilde{\Delta}_{30}(z):=\Delta_{35}(z)/\Delta_5(2z)$.
In the orthogonal language of Sect. 2, this is the automorphic
form $\widetilde{\Delta}_{30}=\widetilde{\Delta}(L_{1,0})_{30}$
introduced in Sect. 2.
Using the lattice $M_{1,I}$, we can rewrite
$\widetilde{\Delta}_{30}(z)$ as follows:

\vskip5pt
$\widetilde{\Delta}_{30}(z)=\Delta_{35}(z)/\Delta_5(2z)=$
$$
=\exp{\left(-2\pi i(\rho_{I},z)\right)}
\prod_{\alpha \in \Delta_+(M_{1,I})^{\re}_+\cup
\Delta(M_{1,0}^\ast )^{\im}_+}
{\left( 1-\exp{ \left(-2\pi i (\alpha ,
z)\right)}\right)^{\mult~\alpha}}
\tag3.10
$$
where
$$
\mult~\alpha = f_2(-2(\alpha, \alpha))-
\cases
f(-(\alpha, \alpha)/2) &\text{if $\alpha \in
M_{1,II}=2M_{1,0}^\ast$,}\\
0                      &\text{if $\alpha \notin
M_{1,II}=2M_{1,0}^\ast$}.
\endcases
\tag3.11
$$
Using these preliminary calculations, we get

\proclaim{Theorem 3.1} The automorphic form $\Delta_{35}(z)$
gives the automorphic correction of the Kac--Moody
algebra $\geg (A_{1,0})$ with the product expansion
\thetag{3.8}, \thetag{3.9}.

The automorphic form
$\widetilde{\Delta}_{30}(z)=\Delta_{35}(z)/\Delta_5(2z)$
gives the automorphic correction of the Kac--Moody
algebra $\geg (A_{1,I})$ with the product expansion \thetag{3.10},
\thetag{3.11}.

The automorphic form $\Delta_5(2z)$
gives the automorphic correction of the Kac--Moody
algebra $\geg (A_{1,II})$ with the product expansion
\thetag{3.7}, \thetag{3.6}.
\endproclaim
\demo{Proof} The proof is similar to the proof of Theorem
2.3 in \cite{GN1}.
For example, consider $\geg (A_{1,0})$. From the
product \thetag{3.8}, one has that
$\Delta_{35}(s_\alpha(z))=-\Delta_{35}(z)$ for any
$\alpha \in \Delta (M_{1,0})^{\re}_+$ because $\mult~\alpha =1$.
It follows that $\Delta_{35}(w(z))=\det (w)\Delta_{35}(z)$,
$w \in W=W^{(2)}(M_{1,0})$.
Since all $\mult~\alpha$ are integral, it then follows from
the product that
$$
\Delta_{35}(z)=\sum_{w\in W}{\det(w)
\left(
-\sum_{ \rho_0 +a \in M_{1,0}^\ast \cap \br_+\M_{1,0}}
{m(a)\exp(-2\pi i (w(\rho_0+a),z))}\right)}
$$
where all $m(a)$ are integral and $m(0)=-1$.
Consider $\rho_0 +a \in M_{1,0}^\ast \cap \br_+\M_{1,0}$.
We have $(\rho_0+a, \delta_i)\le 0$, $i=1,2,3$. If
$(\rho_0+a, \delta_i)=0$, then the corresponding Fourier coefficient
$m(a)=0$, since $\Delta_{35}(z)$ is anti-invariant with respect to
$s_{\delta_i}$. Thus, we can suppose that $(\rho_0+a, \delta_i)<0$,
$i=1,2,3$, considering only non-zero $m(a)$.
By definition of $\rho_0$, we
then get that $(a, \delta_i)\le 0$ since $a\in M_{1,0}^\ast$.
By \thetag{3.1}, then $a \in M_{1,0}^\ast \cap \br_+\M_{1,0}$.
It follows that $a \in \br_{++}\M_{1,0}$ if $a \not=0$.
If $a=0$, we have $m(a)=-1$. Thus $\Delta_5(z)$ has the
form \thetag{3.3}. It follows the statement.
\newline
\qed
\enddemo
\remark{Remark 3.2} In Sects. 1 and 2 we have considered the
Igusa modular  form
$\Delta_{30}(z)=\Delta_{35}(z)/\Delta_5(z)=c\chi_{30}$.
It seems $\Delta_{30}(2z)$ with the product expansion of
Corollary 1.6 gives an automorphic
correction of the Kac--Moody superalgebra
with the generalized Cartan matrix $A_{1,0}$ and with the
set of odd indexes $\tau=\{ 2 \} \subset I=\{1,2,3\}$.
See \cite{K3} and  \cite{R} about these algebras.
We don't consider automorphic corrections of
Kac--Moody superalgebras in this paper.
\endremark

\head
4. The perspective
\endhead

Using methods and results of \cite{N5}, \cite{N8} (see also
general results in \cite{N3}, \cite{N4}), we can prove
the following classification results:

\proclaim{Theorem 4.1} There are exactly 12 elliptic
hyperbolic symmetric generalized Cartan matrices
of rank $3$ which have the Weyl group with a non-compact (i.e.
with an infinite vertex)
fundamental polygon and have a lattice Weyl vector.
They are matrices $A_{1,0}$ -- $A_{3,III}$ below:
$$
\align
A_{1,0}&=
\pmatrix
\hphantom{-}{2}&\hphantom{-}{0}&{-1}\cr
\hphantom{-}{0}&\hphantom{-}{2}&{-2}\cr
{-1}&{-2}&\hphantom{-}{2}\cr
\endpmatrix,\qquad
A_{1,I}=
\pmatrix
\hphantom{-}{2}&{-2}&{-1}\cr
{-2}&\hphantom{-}{2}&{-1}\cr
{-1}&{-1}&\hphantom{-}{2}\cr
\endpmatrix,\\
A_{1,II}&=
\pmatrix
\hphantom{-}{2}&{-2}&{-2}\cr
{-2}&\hphantom{-}{2}&{-2}\cr
{-2}&{-2}&\hphantom{-}{2}\cr
\endpmatrix,\quad\,
A_{1,III}=
\pmatrix
\hphantom{-}{2}&{-2}&{-6}&{-6}&{-2}\cr
{-2}&\hphantom{-}{2}&\hphantom{-}{0}&{-6}&{-7}\cr
{-6}&\hphantom{-}{0}&\hphantom{-}{2}&{-2}&{-6}\cr
{-6}&{-6}&{-2}&\hphantom{-}{2}&{0}\cr
{-2}&{-7}&{-6}&\hphantom{-}{0}&\hphantom{-}{2}\cr
\endpmatrix ,
\endalign
$$
$$
A_{2,0}=
\pmatrix
\hphantom{-}{2}&{-2}&{-2}\cr
{-2}&\hphantom{-}{2}&\hphantom{-}{0}\cr
{-2}&\hphantom{-}{0}&\hphantom{-}{2}\cr
\endpmatrix,
$$
$$
A_{2,I}=
\pmatrix
\hphantom{-}{2}&{-2}&{-4}&\hphantom{-}{0}\cr
{-2}&\hphantom{-}{2}&\hphantom{-}{0}&{-4}\cr
{-4}&\hphantom{-}{0}&\hphantom{-}{2}&{-2}\cr
\hphantom{-}{0}&{-4}&{-2}&\hphantom{-}{2}\cr
\endpmatrix ,
\qquad
A_{2,II}=
\pmatrix
\hphantom{-}{2}&{-2}&{-6}&{-2}\cr
{-2}&\hphantom{-}{2}&{-2}&{-6}\cr
{-6}&{-2}&\hphantom{-}{2}&{-2}\cr
{-2}&{-6}&{-2}&\hphantom{-}{2}\cr
\endpmatrix ,
$$
$$
A_{2,III}=
\pmatrix
\hphantom{-}{2}&{-2}&{-8}&{-16}&{-18}&{-14}&{-8}&\hphantom{-}{0}\cr
{-2}&\hphantom{-}{2}&\hphantom{-}{0}&{-8}&{-14}&{-18}&{-16}&{-8}\cr
{-8}&\hphantom{-}{0}&\hphantom{-}{2}&{-2}&{-8}&{-16}&{-18}&{-14}\cr
{-16}&{-8}&{-2}&\hphantom{-}{2}&\hphantom{-}{0}&{-8}&{-14}&{-18}\cr
{-18}&{-14}&{-8}&\hphantom{-}{0}&\hphantom{-}{2}&{-2}&{-8}&{-16}\cr
{-14}&{-18}&{-16}&{-8}&{-2}&\hphantom{-}{2}&\hphantom{-}{0}&{-8}\cr
{-8}&{-16}&{-18}&{-14}&{-8}&\hphantom{-}{0}&\hphantom{-}{2}&{-2}\cr
\hphantom{-}{0}&{-8}&{-14}&{-18}&{-16}&{-8}&{-2}&\hphantom{-}{2}\cr
\endpmatrix ,
$$
$$
A_{3,0}=
\pmatrix
\hphantom{-}{2}&{-2}&{-2}\cr
{-2}&\hphantom{-}{2}&{-1}\cr
{-2}&{-1}&\hphantom{-}{2}\cr
\endpmatrix ,
\qquad\qquad\ A_{3,I}=
\pmatrix
\hphantom{-}{2}&{-2}&{-5}&{-1}\cr
{-2}&\hphantom{-}{2}&{-1}&{-5}\cr
{-5}&{-1}&\hphantom{-}{2}&{-2}\cr
{-1}&{-5}&{-2}&\hphantom{-}{2}\cr
\endpmatrix ,
$$
$$
A_{3,II}=
\pmatrix
\hphantom{-}{2}&{-2}&{-10}&{-14}&{-10}&{-2}\cr
{-2}&\hphantom{-}{2}&{-2}&{-10}&{-14}&{-10}\cr
{-10}&{-2}&\hphantom{-}{2}&{-2}&{-10}&{-14}\cr
{-14}&{-10}&{-2}&\hphantom{-}{2}&{-2}&{-10}\cr
{-10}&{-14}&{-10}&{-2}&\hphantom{-}{2}&{-2}\cr
{-2}&{-10}&{-14}&{-10}&{-2}&\hphantom{-}{2}\cr
\endpmatrix ,
$$
$$
A_{3,III}=
$$
$$
\pmatrix
\hphantom{-}{2}&{-2}&{-11}&{-25}&{-37}&{-47}&{-50}
&{-46}&{-37}&{-23}&{-11}&{-1}\cr
{-2}&\hphantom{-}{2}&{-1}&{-11}&{-23}&{-37}&
{-46}&{-50}&{-47}&{-37}&{-25}&{-11}\cr
{-11}&{-1}&\hphantom{-}{2}&{-2}&{-11}&{-25}
&{-37}&{-47}&{-50}&{-46}&{-37}&{-23}\cr
{-25}&{-11}&{-2}&\hphantom{-}{2}&{-1}&{-11}
&{-23}&{-37}&{-46}&{-50}&{-47}&{-37}\cr
{-37}&{-23}&{-11}&{-1}&\hphantom{-}{2}&{-2}
&{-11}&{-25}&{-37}&{-47}&{-50}&{-46}\cr
{-47}&{-37}&{-25}&{-11}&{-2}&\hphantom{-}{2}
&{-1}&{-11}&{-23}&{-37}&{-46}&{-50}\cr
{-50}&{-46}&{-37}&{-23}&{-11}&{-1}&\hphantom{-}{2}
&{-2}&{-11}&{-25}&{-37}&{-47}\cr
{-46}&{-50}&{-47}&{-37}&{-25}&{-11}&{-2}
&\hphantom{-}{2}&{-1}&{-11}&{-23}&{-37}\cr
{-37}&{-47}&{-50}&{-46}&{-37}&{-23}&{-11}&{-1}
&\hphantom{-}{2}&{-2}&{-11}&{-25}\cr
{-23}&{-37}&{-46}&{-50}&{-47}&{-37}
&{-25}&{-11}&{-2}&\hphantom{-}{2}&{-1}&{-11}\cr
{-11}&{-25}&{-37}&{-47}&{-50}&{-46}&{-37}
&{-23}&{-11}&{-1}&\hphantom{-}{2}&{-2}\cr
{-1}&{-11}&{-23}&{-37}&{-46}&{-50}
&{-47}&{-37}&{-25}&{-11}&{-2}&\hphantom{-}{2}\cr
\endpmatrix .
$$
\endproclaim

\proclaim{Theorem 4.2} Let $M_{4,0}=U(16)\oplus \langle 2 \rangle$.
Let $\M_{4,0}$ be a fundamental  polygon for
$W=W^{(2)}(M_{4,0})$ in $\La^+(M_{4,0})$ and
$P(\M_{4,0})\subset \Delta^{(2)}(M_{4,0})$ the set of all orthogonal
vectors to faces of $\M_{4,0}$.

The Gram matrix $A_{4,0}=G(P(\M_{4,0}))$ of elements of
$P(\M_{4,0})$ (this matrix is infinite)
is the only parabolic hyperbolic symmetric generalized Cartan
matrix of rank $3$ which has the Weyl group with a
fundamental polygon which has at least one infinite vertex
different from the cusp and has a lattice Weyl vector.
\endproclaim

Using methods developed in \cite{GN1}-- \cite{GN3}
and in this paper, we can find automorphic corrections for all
Lorentzian Kac--Moody algebras $\geg (A)$ where $A$ is
one of 13 generalized Cartan matrices of Theorems 4.1 and 4.2.
There are analogues of automorphic forms $\Delta_5$, $\Delta_{35}$,
$\Delta_{30}$, $\widetilde{\Delta}_{30}$,
$\widetilde{\widetilde{\Delta}}_{30}$ for
all matrices $A_{i,II}$, $A_{i,0}$, $A_{i,I}$, $i=1,\dots,4$,
respectively. The cases $A_{1,II}$ and $A_{2,II}$ had been
considered in \cite{GN1}, and the cases $A_{1,0}$, $A_{1,I}$ have
been considered in this paper. The cases $A_{i,III}$, $i=1,2,3$,
are more delicate.
Using multiplicative Hecke operators applied to these forms,
one can find analogues of modular  forms $F_p(Z)$ which we
consider in Sect. 1, and modular  forms of
Appendix B, and their product and sum expansions.
We hope to publish these results later.

\head
Appendix A.
Multiplicative  Hecke operators and the lifting of Jacobi forms
\endhead
In Appendix  we give a  proof of a  generalization
of  Theorems 1.5 and 1.7.
In particular, this  provides us with a new proof of these theorems
and  explains why the Hecke-Jacobi operator $T^J(p)$ appears
in the formula \thetag{1.12} for the exponents in the infinite
product expansion of the modular form $\dtp$.

For an arbitrary positive integer $t$, we denote by
$$
\gmt=\left\{\pmatrix *&t*&*&*\\
*&*&*&t^{-1}*\\
*&t*&*&*\\*&t*&t*&*
\endpmatrix \in Sp_2(\Bbb Q)\right\},
$$
the paramodular group of type $(1,t)$
where all $*$ are integers.
(The threefold $\gmt\setminus \bh$ is a coarse moduli
space of abelian surfaces with a polarization of type $(1,t)$.)
For $t=1$, we have $\gm=Sp_4(\bz)$.

\definition{Definition} A holomorphic function
$\phi(z_1, z_2): \Bbb H_1\times \bc\to \bc$
is called a {\it Jacobi form} of index $t\in \bn$ and weight $k$
if the function
$\widetilde\phi(Z)=\phi(z_1,z_2)\hbox{\rm exp}(2\pi i\,t z_3)$
on the Siegel upper half-plane $\bh$
is a modular form of weight $k$ with respect to the  parabolic
subgroup $\gi$ (see \thetag{1.7}).
It means that
$\widetilde\phi(Z)$ satisfies the functional equation \thetag{1.1}
for any $M\in \gi$ and that
$\widetilde\phi(Z)$ is holomorphic at infinity,
i.e. it has a Fourier expansion of the type $$
\widetilde\phi(Z)=\sum\Sb n, l\in \bz \\
\vspace{0.5\jot}4nt-l^2
\ge 0\endSb
f(n,l)\,\exp{(2\pi i(n z_1+lz_2+t z_3))}.
$$ \enddefinition
We remark that a Jacobi form of index $t$ satisfies
the functional equation \thetag{1.1} for any
$M\in \Gamma_\infty^{(t)}= \gmt\cap \gi(\bq)$.
We shall use the term Jacobi form for the function
$\widetilde\phi(Z)$ as well.

We call the  function  $\phi(z_1, z_2)$  (or $\widetilde\phi(Z)$)
a  {\it Jacobi cusp form} if we have the strict
inequality $4nt> l^2$ in the Fourier expansion.

The function $\phi(z_1, z_2)$  (or $\widetilde\phi(Z)$)
is called a {\it weak Jacobi form}
if we have the condition $n\ge 0$ instead of  $4nt\ge l^2 $ in  the

Fourier expansion (see \cite{EZ, \S 9}).

The space of  all Jacobi forms
(resp. Jacobi cusp forms) of index $t$ and weight $k$ is denoted by

$M_{k,t}^J$ (resp. $S_{k,t}^J$).

The maximal parabolic subgroup $\gi$ is the semi-direct product
of $SL_2(\bz)$ and the integral  Heisenberg group$$
\gi/\{\pm E_4\}\cong
\left\{\pmatrix a&0&b&0\\
0&1&0&0\\
c&0&d&0\\
0&0&0&1
\endpmatrix \right\}\ltimes
\left\{\pmatrix 1&0&0&\mu\\
\lambda&1&\mu&\kappa\\
 0&0&1&-\lambda\\
0&0&0&1 \endpmatrix
\right\}
$$
where  $\pmatrix a&b\\c&d\endpmatrix\in SL_2(\bz)$
and $\mu, \lambda, \kappa\in \bz$.
Thus, for Jacobi forms, the equation $\thetag{1.1}$is equivalent to

two separate  functional equations:
one for  elements of $\slt$, and another for elements of
the Heisenberg group.

The Hecke algebra $\Cal H(\gi)=\Cal H_\bq(\gi , \hbox{G}\gi)$
will be very useful in our consideration.
As a linear space over $\bq$, this algebra is generated bythe double

cosets
$$
U=\gi M \gi=\sum_i \gi M_i\in \Cal H(\gi)
$$
where $M$ and $M_i$ are  elements of the group of integral
symplectic similitudes $\hbox{G}\gi$ of type $\gi$
(see \thetag{1.2}).
The algebra $\Cal H(\gi)$ is  equipped  with the usual
multiplication
(see \cite{G3}, \cite{G5}   for a general point of view
on such algebras).

This  Hecke algebra acts on the space of  Jacobi
forms of all indices.
Let $V=\sum_i  a_i\, \gi M_i\in \Cal H(\gi)$ be an arbitrary
element and
let $F(Z)$ be any  function which is invariant  with respect
to the  $|_k\,$-action  of the parabolic subgroup   $\gi$,
i.e. $F|_k \gamma= F$ for any $\gamma\in \gi$ (see \thetag{1.1}).
Then we put
$$
F(Z)\to (F|_k\,V)(Z)=
\cases  \sum_i \mu(M_i)^{2k-3} a_i\, (F|_k\,M_i)(Z)
&\text{if } k>0\\
\sum_i a_i\, (F|_k\,M_i)(Z)
&\text{if } k=0
\endcases
\tag{A.1}
$$
where $\mu(M_i)$ denotes the degree of the symplectic similitudes
(see \thetag{1.2}).
(In the definition for, $k>0$, we keep  the same normalizing
factor as for the Hecke operators for $Sp_4(\Bbb Z)$.)

If $\widetilde\phi(Z)=\phi(z_1,z_2)\exp{(2\pi i\, t z_3)}$
is a Jacobi form of weight $k$ and index $t$, then
$(\widetilde\phi|_k\,V)(Z)$ is a sum of Jacobi forms
of possibly different  indices.
To describe this representation, we need to define the concept
of  signature of double cosets  from $\Cal H(\gi)$.

Let
$$
U=\gi \pmatrix *&0&*&*\\
                  *&a&*&*\\
                  *&0&*&*\\
                  0&0&0&d \endpmatrix \gi
\in \Cal H(\gi)  \qquad (a,d\in \bn).
$$
The rational number $\hbox{Sign\,}(U)=ad^{-1}$ is called
the {\it signature} of $U$.
It is a multiplicative invariant of the double cosets
since
$$
\hbox{Sign\,}(\gi M_1\gi\cdot \gi M_2\gi)
=\hbox{Sign\,}(\gi M_1\gi)\cdot \hbox{Sign\,}(\gi M_2\gi).
$$
One can easily prove
(see \cite{G5, Proposition 5.2}).

\proclaim{Lemma A.1}Let $\widetilde\phi(Z)
=\phi(z_1,z_2)\exp{(2\pi  i\, t z_3)}\in M^J_{k,t}$,
and $U\in \Cal H(\gi)$ be a double coset of  signature
$ad^{-1}$.
Then we have $(\widetilde\phi|_k U)(Z)\in M^J_{k,tad^{-1}}$ if
$tad^{-1}\in \bn$ and $(\widetilde\phi|_k  U)(Z)\equiv 0$ otherwise.
\endproclaim

The Hecke algebra $\Cal H(\gi)$ is not commutative and contains  zero

divisors. This algebra has two important properties:

a) $\Cal H(\gi)$ might be viewed as a non-commutative extension  of

the standard
\newline
Hecke algebra $\Cal H (\gm)=\Cal H_\bq(\spi , GSp_4(\bz))$ of the
Siegel modular group;

b) $\Cal H(\gi)$ contains two copies of the Hecke algebra
of the special linear group
$\Cal H(\slt)=\Cal H_\bq(\slt, M_2^+(\bz))$.

These two  algebras are the tensor products of their local components

$$
\Cal H(\slt)=\bigotimes_p \Cal H_p(\slt), \quad \Cal H(\gm)
=\bigotimes_p \Cal H_p(\gm )
$$
which   are the polynomial rings
with generators
$$
T(p)=\slt \hbox{diag}(1,p)\slt,
\quad T(p,p)=\slt\, pE_2\,\slt
$$
and
$$
T_p=\gm \hbox{diag}(1,1,p,p)\gm ,\quad T_{1,p}=\gm
\hbox{diag}(1,p,p^2,p)\gm ,\quad  \Delta_{p}=\gm (pE_4)\gm
$$
respectively.
We remind that
$$
\Cal H_p(\Gamma)
=\{\sum_i a_i\Gamma M_i\in \Cal H_p(\Gamma)\ |\
\mu(M_i)=p^{m_i}\}
$$
where $\Gamma$ is one of the groups $\gmt$, $\gi$ or
$\Gamma_\infty^{(t)}$.

There is a natural  embedding of $\Cal H(\gm)$ into $\Cal H(\gi)$.
If  $V=\sum_i a_i\,\gm M_i\in \Cal H(\gm)$,  then according to the

elementary divisors theorem   one can represent $V$ in the form
$V=\sum_i a_i\,\gm M'_i$, where $M'_i \in \hbox{G}\gi(\bz)$.
It is  easy to see that the map
$$
Im: U=\sum_i a_i\,\gm M_i\to \sum_i
a_i\,\gi M'_i
$$
is an embedding of $\Cal H(\gm)$ into $\Cal
H(\gi)$ (see \cite{G3}, \cite{G5} for  more general constructions).

We shall  identify  $\Cal H(\gm)$ with its image in $\Cal H(\gi)$.

The ring $\Cal H (\gi)$  contains  two subrings  isomorphic to

$\Cal H (\slt)$:
$$
\CD
@. \Cal H(SL_2(\bz)) \\
@.  @Vj_-Vj_+V\\
\Cal H(\gm) @>Im>> \Cal H(\gi).
\endCD
$$
It is enough  to define the embeddings $j_{\pm}$
for the generators
$T(p)$ and $T(p,p)$ whose images we denote by
$T_{\pm}(p)=j_{\pm}(T(p))$ and
$\Lambda_{\pm}(p)= j_{\pm}(T(p,p))$.
By definition, we have
$$
\align
T_{-}(p)&=\gi\text{diag}(1,p,p,1)\gi,\quad
\Lambda_-(p)=\gi\text{diag}(p,p^2,p,1)\gi,\\
T_{+}(p)&=\gi\text{diag}(1,1,p,p)\gi,\quad
\Lambda_+(p)=\gi\text{diag}(p,1,p,p^2)\gi.
\endalign
$$
These elements are connected with the standard generators of
the symplectic Hecke algebra $\Cal H_p(\gm)$
$$
T(p)=T_{-}(p)+T_{+}(p),\quad
T_{1,p}=\Lambda_{-}(p)+\Lambda_{+}(p)
+T_0(p)+\nabla_p-\Delta_p
\tag{A.2}
$$
where
$$
T_0(p)=\gi\text{diag}(1,p,p^2,p)\gi,\quad
\nabla_p=\sum_{r\in \bz/p\bz} \gi
\pmatrix
p&0&0&0\\
0&p&0&r\\
0&0&p&0\\
0&0&0&p
\endpmatrix
$$
and $\Delta_p=\gi (pE_4) \gi$.
(The last element is a central element of $\Cal H (\gi)$.)

For any
$V=\sum_i a_i\,\gi M_i\in \Cal H(\gi)$, we denote by
$\hbox{Sign}_s(V)$ or $V^{(s)}$
its  homogeneous part of signature $s$
$$
V^{(s)}=\hbox{Sign}_s(V)=\sum_{i,\ {Sign}(M_i)=s} a_i\, \gi M_i.
$$

The commutativity of the elements $T_\pm(p)$, $\Lambda_\pm(p)$,
$T_0(p)$ and $\nabla_p$,
for different primes, follows from the  lemma
(see \cite{G2, Satz 5.1}).

\proclaim{Lemma A.2} Let $U=\sum_{a} U^{(a)}$ and
$V=\sum_{b} V^{(b)}$
be the decompositions in $\Cal H(\gi)$
of  elements $U, V\in \Cal H(\gm)$ in the sum of homogeneous
components with the same signature.
Then the  identity
$$
\sum_{ab=r} U^{(a)}\cdot V^{(b)}=
\sum_{ab=r} V^{(b)}\cdot U^{(a)}
$$
is valid.
\endproclaim

We shall see that the formula of Theorem 1.7
for $f_p(N)$ and the formula  \thetag{1.12} follow
from  commutator relations in the non-commutative ring
$\Cal H(\gi)$, which reflect some properties of  local
$L$-functions of $\spi$ and $\slt$.

We remind that the Hecke element  $T(m)\in \Cal H(\slt)$
is equal to the sum of all distinct double cosets
with determinant $m$
$$
T(m)=\sum_{a|b,\, ab=m}
\slt \pmatrix a&0\\0&b\endpmatrix  \slt.
$$
For a prime $p$, we
introduce the $p$-local Hecke polynomial
$$
Q_p(X)=1-T(p)X+p\,T(p,p)X^2
$$
which is the denominator of the local
zeta-function of the group $\slt$
$$
\sum_{d\ge 0}T(p^d)\, X^d=Q_p(X)^{-1}.
\tag{A.3}
$$
A similar to  $Q_p(X)$ polynomial  for the symplectic group
is the  polynomial over $\Cal H_p(\gm)$
$$
S_p(X)=
1-T_pX+p\bigl(T_{1,p}+(p^2+1)\Delta_p\bigl)X^2
-p^3\Delta_pT_pX^3
+p^6\Delta_p^2X^4
$$
which  is the denominator of the Hecke series of type \thetag{A.3}
for $Sp_4(\bz)$.
This polynomial is the local factor of the so-called
Spin-$L$-function of $Sp_4(\bz)$.

The next identity (see \cite{G5}) shows a connection between
the polynomial $S_p(X)$ and  the $\pm$-embeddings
of $Q_p(X)$  into the parabolic extension $\Cal H(\gi)$
of $\Cal H(\gm)$. It is
$$
S_p(X)=Q_{p,-}(X)\bigl(1-K_p\, X^2\bigr)Q_{p,+}(X)
\tag{A.4}
$$
where
$$
Q_{p,\pm}(X)=1-T_\pm (p)X+p\Lambda_{\pm}(p)X^2,
\qquad K_p= p^2\Delta_p-p\nabla_p.
$$
(The proof of \thetag{A.4} and some more general identities
of this type see in \cite{G3}, \cite{G5}.)

The ``minus''-embedding is used in the construction of the
arithmetical
lifting of the Jacobi forms.
By  definition of the ``minus''-embedding, we have
$$
T_{-}(m)=\sum\Sb ad=m\\
\vspace{0.5\jot} b\,mod\,d\endSb
\gi\pmatrix a&0&b&0\\0&m&0&0\\0&0&d&0\\0&0&0&1\endpmatrix.
$$
Hence  by Lemma A.1
$$
\bigl(\widetilde \phi|_k\, \,T_{-}(m)\bigr)(Z)=
m^{2k-3} \sum\Sb ad=m\\ \vspace{0.5\jot} b\,mod\,d\endSb d^{-k}
\phi (\frac{az_1+b}d,\, az_2)
\exp{(2\pi i\, mt z_3)} \in M^J_{k,mt}
$$
for any
$\widetilde\phi(Z)=\phi(z_1,z_2)\exp{(2\pi i\, t z_3)}
\in M^J_{k,t}$ ($k>0$).
In \cite{G1, Theorem 3}
it was proved that for any Jacobi cusp form
$\widetilde\phi(Z)\in S^J_{k,t}$ ($k>0$) the function
$$
\hbox{Lift\,}(\widetilde\phi)(Z)=
\sum_{m=1}^{\infty}\,m^{2-k}\,
\bigl(\widetilde\phi\,|_{k}\,T_{-}(m)\bigr)(Z)
\tag{A.5}
$$
is a cusp form of weight $k$ with respect
to the paramodular group $\gmt$.

Another application of the same method
was found in \cite{GN1} where we constructed
modular  forms with respect to the paramodular group
$\gmt$ of type
$$
G(Z)=\widetilde\psi(Z)
\exp{\biggl(-\sum_{m=1}^{\infty}\,m^{-1}\,
\bigl(\widetilde\phi\,|_{0}\,T_{-}(m)\bigr)(Z)\biggr)}
\tag{A.6}
$$
where $\widetilde\psi(Z)$ is a Jacobi cusp form
of integral or  half integral index and
$\widetilde\phi(Z)$ is a weak Jacobi form
of weight $0$ and index $t$.
The identity of type \thetag{A.6} implies the product formula
\thetag{0.2} of the modular form $\Delta_5(Z)$
(see \cite{GN1, Theorem 4.1}).

For $t=1$ the lifting \thetag{A.5} coincides with
the Maass lifting and it commutes with
action of Hecke operators from $\Cal H(\gm)$.
A proof of this fact in  \cite{EZ, \S 6}
is based on the exact  formulae for action of
Hecke operators on Fourier coefficients of Siegel
modular forms. Another  proof using  parabolic Hecke operators
from $\Cal H(\gi)$, one can find in \cite{G4}, where
the commutativity of the lifting and Hecke operators
was proved for the Hermitian modular group $SU(2,2)$.
In \cite{G4} we used that  the lifting is a modular form.

Bellow we prove  the commutativity of the lifting of Jacobi
forms with action of Hecke operator as a corollary of
some  commutator identities in the non-commutative extension
$\Cal H(\gi)$. Such relations provide a purely algebraic
proof of the result mentioned above for any paramodular group
$\gmt$  and  an exact formula for
the action of arbitrary multiplicative Hecke operator of type
\thetag{1.3} on modular forms of type \thetag{A.6}.

The next proposition will  play the crucial role in the proof.

\proclaim{Proposition A.3}1. For any element $V\in \Cal H (\gm)$,
the formal power series
$$
R_p(V,X)=Q_{p,-}(X)^{-1}\,V\,Q_{p,-}(X)=
\biggl(\sum_{d\ge 0}T_{-}(p^d)X^d\biggr)
\, V \,Q_{p,-}(X)
$$
is a polynomial in $X$ over $\Cal H(\gi)$.

2. Let $V\in \Cal H_p(\gm)$ and let
$F_t(Z)=\phi(z_1,z_2)\exp{(2\pi i\, tz_3)}$ be
a $|_k$-invariant function with respect to $\gi$  where $(t,p)=1$.
Then
$$
F_t|_k\, R_p(V,\,p^{2-k})=
\sum_{d\ge 0} F_t|_k\,
\hbox{\rm Sign}_1\bigl(T_{-}(p^d) V\bigr)\, p^{(2-k)d}
$$
where $\hbox{\rm Sign}_1(U)$ denote
the part of  $U\in \Cal H_p(\gi)$ with signature $1$ (see Lemma
A.2).
\endproclaim
\demo{Proof}1. $\Cal H(\gm)$ is commutative, hence we have
$$
S_p(X)^{-1}\,V\, S_p(X)=V
$$
for any $V\in \Cal H(\gm)$. Using \thetag{A.4}, we get
$$
Q_{p,-}(X)^{-1}\,V\,Q_{p,-}(X)=
\bigl(1-K_p\, X^2\bigr)Q_{p,+}(X)\,V\,
Q_{p,+}(X)^{-1}\bigl(1-K_p\, X^2\bigr)^{-1}.
\tag{A.7}
$$
By definition, we have
$$
\hbox{Sign}(T_{-}(p^d))=p^d, \quad
\hbox{Sign}(T_{+}(p^d))=p^{-d},\quad
\hbox{Sign}(K_p)=1.
$$
Comparing the signature of the both sides of \thetag{A.7},
we see that the formal power series $Q_{p,-}(X)^{-1}\,V\,Q_{p,-}(X)$
is a polynomial, because
$$
\hbox{min}(\hbox{Sign}(V))\le
\hbox{Sign\,}\biggl(Q_{p,-}(X)^{-1}\,V\,Q_{p,-}(X)\biggr)
\le \hbox{max}(\hbox{Sign}(V))
$$
where $\hbox{min}(\hbox{Sign}(V))$
(resp. \hbox{max}(\hbox{Sign}(V))) is the minimal
(resp. maximal) signature of the components of $V$.

2.  According to Lemma A.1 and to the condition $(t,p)=1$,
we have $F_t|_k W\equiv 0$ if $\hbox{Sign}(W)=p^{-n}$
with $n>0$.
Hence all elements with signature of type $p^{-n}$ $(n\in \bn)$
disappear after $|_k$-action of $V$ on $F_t$.
Elements with signature  $p^{n}$ multiply
the index $t$ of the function $F_t$.
It is easy to see that
$$
F_d|_k K_p=F_d|_k(p^2\Delta_p-p\nabla_p)=
\cases p^{2k-4}F_d &\text{if } (d,p)=1,\\
 0&\text{if } (d,p)=p.
\endcases
$$
Let us write \thetag{A.7} in the  form
$$
\align
R_p(V,X)&=Q_{p,-}(X)^{-1}\,V\,Q_{p,-}(X)\\
{}&=\bigl(1-K_p\, X^2\bigr)Q_{p,+}(X)\,V\,
Q_{p,+}(X)^{-1}+
Q_{p,-}(X)^{-1}\,V\,Q_{p,-}(X) K_p\, X^2.
\endalign
$$
This identity implies the second statement  since
$F_t|_k (1-K_p X^2)\equiv 0$ for $X=p^{2-k}$
and $F_t|_k\,\bigl(\hbox{Sign}_{p^n}(R_p(V,X)) K_p\bigr)=0$
for $n\ge 1$.
\newline\qed
\enddemo
Let us introduce a  map $\Cal J$ from the Hecke algebra
$\Cal H(\gm)$ into the linear space
$$
\Cal H^{(0)}(\gi)=
\{V=\sum_i\, a_i\, \gi M_i\in \Cal H(\gi)\,|\, \hbox{Sign}(M_i)=1\}.
$$
For any $V=\sum_a V^{(a)}\in \Cal H(\gm)$ (see Lemma A.2) we put
$$
\Cal J_k(V):=\sum_{m\ge 1} m^{2-k}\, T_-(m)V^{(m^{-1})}
\in \Cal H^{(0)}(\gi).
\tag{A.8}
$$
The commutativity between the lifting and Hecke operators
follows from the global variant of the second statement
of Proposition A.3.
\proclaim{Corollary A.4}Let
$\widetilde\phi(Z)=\phi(z_1,z_2)\exp{(2\pi i\, z_3)}\in S^J_{k,1}$
be a Jacobi cusp form of weight $k$ and index $1$ and let
$V\in \Cal H(\gm)$. Then we have
$$
\hbox{\rm Lift\,}(\widetilde\phi)|_k V =
\hbox{\rm Lift\,}(\widetilde\phi|_k \Cal J_k(V)).
$$
\endproclaim
\demo{Proof}Without loss of generality we suppose that
$V$ is a double coset.
Let
$$
V=\gm \hbox{diag}(a,b,c,d)\gm =\prod_{p|v} V_p,
\qquad  (v=ac=bd=\mu(V))
$$
be the  decomposition of $V\in \Cal H(\gm)$
as the product of its local components
$V_p\in \Cal H_p(\gm)$.
The element $T_-(p^d)$ is the component of signature
$p^d$ of the Hecke element $T(p^d)\in \Cal H_p(\gm)$.
Therefore, by Lemma A.2,
$$
\Cal J_k(V)=
\bigl(\prod_{p|\mu(V)} Q_{p,-}(p^{2-k})^{-1}\bigr)\, V\,
\bigl(\prod_{p|\mu(V)} Q_{p,-}(p^{2-k})\bigr)=
\prod_{p|\mu(V)} R_p(V_p,\, p^{2-k}).
$$
Due to the rationality  of the Hecke series of
$\Cal H(\slt)$ (see \thetag{A.3}),
we can rewrite the definition of the lifting
using a formal  infinite product (see \cite{G1})
$$
\hbox{\rm Lift\,}(\widetilde\phi)=
\widetilde\phi|_k \sum_{m=1}^{\infty}\,m^{2-k}\,T_{-}(m)=
\widetilde\phi|_k \prod_p Q_{p,-}(p^{2-k})^{-1}.
$$
In view of Proposition A.3 and Lemma A.2 we get
the formula of Corollary
$$
\hbox{\rm Lift\,}(\widetilde\phi)|_V\kern-1pt =
\kern-1pt \biggl(\widetilde\phi|_k
\bigl(\prod_{p|\mu(V)} Q_{p,-}(p^{2-k})^{-1} \ V
\prod_{p|\mu(V)} Q_{p,-}(p^{2-k})\bigr)\biggr)
\kern-1pt\bigm|_k \kern-1pt
\prod_p Q_{p,-}(p^{2-k})^{-1}.
$$
\newline\qed
\enddemo
\example{Example A.5}(See \cite{EZ, Theorem 6.3}.)
One can  calculate the images
of the generators of the Hecke algebra $\Cal H(\gm)$
under $\Cal J_k$ using \thetag{A.2}.
One can easily check that
$$
\align
T_{-}(p)T_{+}(p)&=pT_0(p)+(p^3+p^2)\Delta_p,
\tag{A.9}\\
T_{-}(p^2)\Lambda_{+}(p)&=p^3T_0(p)\Delta_p+p^4\Delta_p^2.
\tag{A.10}
\endalign
$$
Hence
$$
\align
\widetilde\phi|_k\Cal J_k(T_p)&=
\widetilde\phi|_k (p^{3-k}T_0(p)+p^{k-1}+p^{k-2}),\\
\widetilde\phi|_k \Cal J_k(T_{1,p})&=
\widetilde\phi|_k \bigl((p+1)T_0(p)+p^{2k-6}(p^2-1)\bigr).
\endalign
$$
The element $T_0(p)\in \Cal H(\gi)$ has the  expansion
in the sum of left cosets
$$
\aligned
\gi\hbox{diag}(1,p,p^2,p)\gi&=\kern-2pt
\sum_{a\,mod\,p}\kern-2pt\gi
{\pmatrix
p^2&0&0&0\\
-pa&p&0&0\\
0&0&1&a\\
0&0&0&p
\endpmatrix}\kern-1pt
+\kern-1pt\sum_{a\,mod\,p}\kern-2pt\gi
{\pmatrix
1&0&c&a\\
0&p&pa&0\\
0&0&p^2&0\\
0&0&0&p\endpmatrix}\\
{}&+
\sum_{b\not\equiv 0\,mod\,p}\
\sum_{a\,mod\,p}
\gi
{\pmatrix
p&0&b&ab\\
0&p&ab&a^2b\\
0&0&p&0\\
0&0&0&p
\endpmatrix}.
\endaligned
$$
Using this system of representatives, one can see that
the operator $|_k T_0(p)$ is
connected with  Hecke-Jacobi operator $T^J(p)$
of Eichler-Zagier
$$
(\phi|_{k, t} T^J(p)) (z_1,z_2)=
p^{k-4}\sum
\Sb M\in SL_2(\bz)\setminus
M_2^+(\bz)\\
\hbox{{\eightpoint det}} (M)=p^2\\
\gcd (M)=1\,or\,p^2\endSb
\ \sum_{X\in \bz^2/p \bz^2}\ \phi|_{k,t}\, M\ |\, X
$$
(see \cite{EZ, \S3}) by the formula
$$
\widetilde\phi|_k T_0(p)=p^{k-3}\widetilde\phi|_k T^J(p).
$$
For the  action of $T_{0}(p)$ on  a Jacobi form
$\widetilde\phi(Z)=\phi(z_1,z_2)\exp{(2\pi i\,z_3)}$
of index $1$, we have the formula
$$
\widetilde\phi|_k T_0(p)(Z)=
\sum_{n,l}f_p(4n-l^2)\exp{(2\pi i(nz_1+lz_2+tz_3))}
$$
with
$$
f_p(N)=
\cases p^{k-3}\bigl(f(p^2N)+p^{k-2}f(N)+p^{2k-3}f(\frac{N}p)\bigr)
&\text{if } k>0,\\
p^3f(p^2N)+pf(N)+f(\frac{N}p)
&\text{if } k=0
\endcases
\tag{A.11}
$$
where we use the notation \thetag{1.8}--\thetag{1.9}.
We remark that  Fourier coefficients $f(n,l)$ of an arbitrary
Jacobi form of index $1$ depend only on the norm
$4n-l^2$.
\endexample
The restriction $t=1$ in Corollary A.4 is not principal.
The lifting \thetag{A.5} is a modular form with respect
to the paramodular group $\gmt$. The Jacobi form
$\widetilde{\phi}_t(Z)=
\widetilde{\phi}(z_1,z_2)\exp{(2\pi i\, tz_3)}$
of weight $k$ and index $t$ is invariant with respect to the
parabolic subgroup
$\Gamma_\infty^{(t)}= \gmt\cap \gi(\bq)$.
Let us consider the Hecke algebra
$\Cal H(\Gamma_\infty^{(t)})$ containing
$\pm$-embeddings of $\Cal H(\slt)$
(one may use the same definition as before) and
the Hecke algebra of all elements having a good reduction
for  the paramodular group $\gmt$
$$
\Cal H_*(\gmt)=\bigotimes_{(p,t)=1} \Cal H_p(\gmt)
\cong \bigotimes_{(p,t)=1} \Cal H_p(\gm).
$$
For a prime $p$, which does not divide $t$, the new Hecke
algebras are isomorphic to the Hecke algebras
connected with the symplectic group
$$
\Cal H_p(\gmt)\cong \Cal H_p(\gm),
\qquad
\Cal H_p(\Gamma_\infty^{(t)})\cong \Cal H_p(\gi).
$$
Thus we can define the same morphism $\Cal J$ for
any element $V=\sum_a V^{(a)}\in \Cal H_*(\gmt)$
$$
\Cal J^{(t)}_k(V)=
\sum_{m\ge 1} m^{2-k}\, T_-(m)V^{(m^{-1})}\in
\Cal H^{(0)}(\Gamma_\infty^{(t)}).
\tag{A.8-t}
$$
Then we have
\proclaim{Corollary A.4-t}Let
$\widetilde\phi_t(Z)=\phi(z_1,z_2)\exp{(2\pi i\, t z_3)}\in
S^J_{k,t}$
be a Jacobi cusp form of weight $k$ and index $t$,
and $V\in \Cal  H_*(\gmt)$ be a Hecke element with good
reduction. Then
$$
\bigl(\hbox{\rm Lift\,}(\widetilde\phi_t)\bigr)|_k V =
\hbox{\rm Lift\,}\bigl(\widetilde\phi_t|_k \Cal J^{(t)}_k(V)\bigr)
$$
\endproclaim
(See \cite{G2} where the local factors of the Spin-$L$-function
of the lifting were calculated for  prime $p$
with an  arbitrary reduction.)

Now we consider an application of Proposition A.3 to the modular
forms of type \thetag{A.6}. Firstly, we have a Jacobi analogues
of Lemma 1.1.
\proclaim{Lemma A.6}Let
$\widetilde{\phi}(z_1,z_2)\exp{(2\pi i\, tz_3)}\in M^J_{k,t}$
be a  Jacobi form of weight $k$  and index $t$, and
$V=\Gamma_\infty^{(t)}M\Gamma_\infty^{(t)}
=\sum_i \Gamma_\infty^{(t)}M_i \in \Cal H(\Gamma_\infty^{(t)})$
be a double coset of the parabolic Hecke ring.
Then the function
$$
[\widetilde{\phi}(Z)]_V:=\prod_i(\widetilde{\phi}|_k M_i)(Z)
$$
is a Jacobi form of weight $k\nu$  and index
$(t\nu)\hbox{\rm Sign(V)}$ where $\nu$
is the number of the left cosets  in $V$.
\endproclaim

\proclaim{Theorem  A.7} Let
$$
G(Z)=\widetilde\psi(Z)
\exp{\biggl(-\sum_{m=1}^{\infty}\,
{\widetilde\phi }\,|_{0}
\,m^{-1}T_{-}(m)(Z)\biggr)}
$$
be  a modular form (see $\thetag{A.6}$)
of weight $k$ with respect to the paramodular
group $\gmt$ where  $\widetilde\psi(Z)$ is a Jacobi form
of weight $k$ and $\widetilde\phi(Z)$  is a Jacobi form of
weight $0$ and index $t$.
Let $V\in \Cal H_*(\gmt)$ be a Hecke element.
Then the  formula
$$
[G(Z)]_V=[\widetilde\psi(Z)]_V
\exp{\biggl(-\sum_{m=1}^{\infty}\,
\bigl({\widetilde\phi_{0,t}| \Cal J_0^{(t)}(V)}\bigr)\,|_{0}
\,m^{-1}T_{-}(m)(Z)\biggr)}
$$
holds where
$$
\Cal J^{(t)}_0(V)=
\sum_{m\ge 1} m^{-1}\, T_-(m)V^{(m^{-1})}\in
\Cal H^{(0)}(\Gamma_\infty^{(t)})
$$
and $[G(Z)]_V$ is defined in \thetag{1.3}.
\endproclaim
\remark{Remark}In the case of Jacobi forms of weight $0$, we
do not have  a normalizing factor in the definition
of the representation $\phi_{0,t}|_0 V$ (see \thetag{A.1}).
This explains  different definition of the operator $\Cal J_0$.
\endremark
\demo{Proof}Let us denote the function under exponent in $G(Z)$
by $F(Z)$.
Each  multiplicative Hecke  operator $[\dots]_V$ acting
on  $\exp{(F(Z))}$ turns into a  usual Hecke operator
on the function $F(Z)$ under  the exponent:
$$
[\exp{(F(Z))}]_V=\exp{(F|_0 V(Z))}.
$$
We remark that the element $V$ to the right is considered
as an element of the Hecke ring $\Cal H(\Gamma_\infty^{(t)})$.
The function $F(Z)$ is  of  the form considered in
Corollary A.4
$$
F(Z)=\widetilde{\phi}|_0 \prod_p Q_{p,-}(p^{-1})^{-1}.
$$
Using Proposition A.3, we get the formula of Theorem A.7.
\newline\qed
\enddemo
\example{Example A.8} Let us calculate $[\Delta_5(Z)]_V$
for $V=T_p$ and $V=T_{1,p}$.
Using \thetag{A.9}--\thetag{A.10}, we have
$$
\align
\widetilde{\phi}|_0 \Cal J_0(T_p)&=
\widetilde{\phi}|_0 \bigl(T_0(p)+p^2+p\bigr)\\
\widetilde{\phi}|_0 \Cal J_0(T_{1,p})&=
\widetilde{\phi}|_0 \bigl((p+1)T_0(p)+p^2-1\bigr).
\endalign
$$
In view of \thetag{A.11}, we get the formula \thetag{1.12}
for the exponents $f_p(N)$ in the infinite product representation
of $\dtp$.
\endexample

\head
Appendix B.
Siegel modular forms with divisors $H_D$.
\endhead

In the remark after Corollary 1.4  we have
mentioned a result of van der Geer.
He  proved in \cite{vdG2}
that there exists a Siegel modular form of weight $-60H(2,D)$
with the divisor
$$
G_D=\bigcup_{e^2|D}H_{D/e^2},
$$
where $H(2,D)$ is the so-called H. Cohen number
and $H_{D/e^2}$ is the Humbert surface of discriminant $D/e^2$.
The proof in \cite{vdG2} is based on the calculation
of degree of the divisor $G_D$ in a smooth compactification of
the Siegel modular threefold
and does not provide an  exact  construction of Siegel
modular forms with known divisors.
Like in Sect. 1 and Appendix, repeating exponential
Hecke operators and considering quotients,
one can easily construct,
starting from the form $\Delta_5(Z)$, Siegel modular forms
with divisors $H_{D^2}$ for any $D$, and their infinite
product expansions.
In this appendix we describe an infinite product construction
of modular forms  with  the divisors  $H_D$
where $D$ is not a perfect square.

An example of the Siegel modular form of weight $24$
with the divisor $H_5$  was found  in \cite{vdG2}
in terms of the Igusa's generators ofthe graded  ring  of Siegel
modular forms:
$$
G^{(5)}=(\chi_{12}-2^{-12}3^{-3}(E_6^2+E_4^3))^2
-E_4(2\cdot 3^{-1}\chi_{10}-2^{11}3^{-6}E_6 E_4)^2.
$$
We find below an  infinite product expansion for this function.
\medskip

Let $j(\tau)$ be the $\slt$-modular invariant function of weight $0$
$$
\align
j(\tau)&=\frac{E_4(\tau)^3}{\Delta_{12}(\tau)}=
\frac{(1+240q+2160q^2+\dots)^3}{q-24q^2+253q^3+\dots}\\
\vspace{2\jot}
{}&=q^{-1}+744+196884q+21493760q^2+\dots\endalign
$$
where $E_4(\tau)$ is the Eisenstein series of weight $4$ and
$\Delta_{12}(\tau)$ is the cusp form of weight $12$
for $\slt$.

The product  $j(z_1)\phi_{0,1}(z_1,z_2)$ (see \thetag{1.8})
is a meromorphic Jacobi form of weight $0$ and index $1$.
This function is not a weak Jacobi form, since it hasthe  Fourier
expansion
$$
\multline
j(z_1)\phi_{0,1}(z_1,z_2)=q^{-1}(r^{-1}+10+r)+\\
+(10r^{-2}+680 r^{-1}+7548+680r+10r^2)+q(\dots).
\endmultline
$$
Let us define another Jacobi form of the same type
$$
\phi^{(5)}_{0,1}(z_1,z_2)=
j(z_1)\phi_{0,1}(z_1,z_2)-10\bigl(\phi_{0,1}|T_0(2)\bigr)(z_1,z_2)
-680\phi_{0,1}(z_1,z_2)
$$
where $T_0(2)$ is the Hecke-Jacobi operator (see \thetag{A.11}).
This Jacobi form has the Fourier expansion of the type
$$
\phi^{(5)}_{0,1}(z_1,z_2)=
\sum\Sb n\ge -1\\\vspace{0.5\jot}
l\in \Bbb Z\endSb
 g(n,l)\,
\exp{(2\pi i \,(nz_1+lz_2))}
=q^{-1}(r^{-1}+r)+48+q(\dots).
$$
Like  in \S 1 we use the Fourier coefficients of the last Jacobi
form  to define an infinite product
which is  a Siegel modular form.
\proclaim{Theorem B.1}The infinite product
$$
F^{(5)}(Z)=r^{-1}(qr-s)(q-rs)\prod\Sb n,\,l,\,m\in \bz\\
\vspace{0.5\jot}
(n,l,m)>0\,\endSb
\bigl(1-q^n r^l s^m\bigr)^{g(nm,l)},
$$
where $(n,l,m)>0$ means that $n\ge 0$, $m\ge 0$, $n+m>0$
and  $l$ is an arbitrary integer,
defines a Siegel modular form of weight $24$. The divisor
of $F^{(5)}(Z)$ is the Humbert surface $H_5$ of discriminant $5$.
\endproclaim
\demo{Proof}One can easily proof this theorem using
the method we used in \cite{GN1, Theorem 4.1}
(see also \cite{B5, Theorem 5.1}).
We remark only the main steps of the proof.

i)  The given  infinite product absolutely converges for all
$\hbox{Im}(Z)>C$ for a large $C>0$ and can be continued
to a multi-valued analytic  function on $\bh$.

ii) It is invariant with respect to the changing
of variables $q\to s$, $s\to q$
(equivalently $z_1\to z_3$, $z_3\to z_1$).

iii) The function $F^{(5)}(Z)$ has a representation
in terms of the Hecke operators $T_-(m)$
$$
F^{(5)}(Z)=\Delta_{12}(z_1)^2
\exp{\biggl(-\sum_{m=1}^{\infty}\,m^{-1}\,
\bigl(\widetilde\phi_{0,1}^{(5)}\,|_{0}\,T_{-}(m)\bigr)(Z)\biggr)}.
$$
Hence the product is a $\gi$-invariant function of weight $24$. The

parabolic subgroup $\gi$ and the transformation

$z_1\to z_3$, $z_3\to z_1$ generate the Siegel modular group.

iv) Only the factors  with negative square $l^2-4nm<0$ of the
infinite product
have zeros on $\bh$.
There factors exist only  with square $l^2-4nm=-5$, and they define

the Humbert surface $H_5$.
\newline
\qed
\enddemo
We remark that the construction given above is general. It gives us
a modular form with the irreducible divisor $H_D$ for any $D$.
For   $D=4d+1$  one can  consider the Jacobi form
$$
j(z_1)^d\cdot
\phi_{0,1}(z_1,z_2)=q^{-d}(r^{-1}+10+r)+q^{-d+1}(\dots).
$$
For $D=4d$  one can start with
the Jacobi form of weight $0$
$$
\bigr(\phi_{0,1}|_0 T_0(2)\bigl)(z_1,z_2)-2\phi_{0,1}(z_1,z_2)=
q^{-1}+(r^2+r^{-2}+70)+q(\dots)
$$
which defines the infinite product expansion  of  $\Delta_{35}(Z)$.

We hope to present  this construction  in more details somewhere.

\enddocument
\end